%% file: arXiv_preprint.tex
\documentclass[11pt]{article}
\usepackage{calc}
\usepackage[a4paper,width=19cm,height=22.5cm]{geometry}
\usepackage[ngerman,english,spanish]{babel}
\usepackage{etex}
\usepackage[T1]{fontenc}
\usepackage[utf8]{inputenc}
\usepackage{ellipsis,ragged2e}
\usepackage{microtype}
\usepackage{amsmath}
\usepackage{amssymb}
\usepackage{amsthm}
\usepackage{color}
\definecolor{rasp}{rgb}{.89,.04,.36}      % rot fuer mid

\usepackage{times}
\usepackage{pictexwd}
\usepackage[pdftex]{graphicx}
\usepackage{graphicx}
\usepackage{caption}
\usepackage{subcaption}
\usepackage{pst-pdf,pstricks-add}
\usepackage{hyperref}
\usepackage{tikz}
\usetikzlibrary{plotmarks}
\usepackage{pgfplots}
\usepackage{dsfont}
\usepackage{kpfonts}
\usepackage{listings}
\usepackage{cleveref}
\usepackage{cite}

\usepackage{mathrsfs,graphicx,color,float, indentfirst,textcomp}
\usepackage{setspace}
\usepackage{latexsym,amssymb,lscape,amsmath,amsthm,amsfonts,amssymb,cite,enumerate}
\usepackage{lscape}

\usepackage{algorithm} %//format of the algorithm
\usepackage{algorithmic} %//format of the algorithm
\usepackage{multirow} %//multirow for format of table
\usepackage{pifont}
\usepackage{doi}

\usepackage{xfrac}
\usepackage{indentfirst}
\usepackage{algorithm} 
\usepackage{algorithmic}
\usepackage{multirow} 
\usepackage{xcolor}
\usepackage{geometry}
\usepackage{booktabs}

\usepackage{tikz}
\usepackage{pgfplots}
\pgfplotsset{compat=newest}

\newcommand\dt{\partial_t}

\definecolor{red}{rgb}{1.0,0.0,0.0}
\definecolor{middlegray}{rgb}{0.5,0.5,0.5}
\definecolor{lightgray}{rgb}{0.8,0.8,0.8}
\definecolor{fhgorange}{rgb}{0.92,0.42,0.04}
\definecolor{fhggreen}{rgb}{0.42,0.61,0.25}
\definecolor{fhgblue}{rgb}{0.0,0.43,0.57}
\definecolor{fhggray}{rgb}{0.66,0.69,0.69}
\definecolor{yac}{rgb}{0.6,0.6,0.1}
\definecolor{luh-dark-blue}{rgb}{0.0, 0.313, 0.608}
\definecolor{luh-light-blue}{rgb}{0.6, 0.725, 0.847}
\definecolor{luh-gray}{rgb}{0.8, 0.863, 0.922}
\definecolor{luh-green}{rgb}{0.784, 0.827, 0.09}

\title{ Hyperbolic Discretization via Riemann Invariants}
\date{May 25, 2019}
\author{Sara Grundel \footnote{grundel@mpi-magdeburg.mpg.de}\\ 
{\small\it Max-Planck-Institut f\"{u}r Dynamik komplexer technischer Systeme} \\
	{\small\it Sandtorstr. 1, 39106 Magdeburg, Germany
}\and  
Michael Herty \footnote{herty@igpm.rwth-aachen.de} \\
	{\small\it Institut f\"{u}r Geometrie und Praktische Mathematik (IGPM)} \\
	{\small\it RWTH Aachen University} \\
	{\small\it Templergraben 55, 52062 Aachen, Germany
}
}

\usepackage{tikz} 
\usetikzlibrary{matrix,positioning,decorations.pathreplacing}

\begin{document}
\selectlanguage{english}
\pagestyle{plain}
%\twocolumn
\date{\today}
\maketitle
% computational and mathematical models
\begin{abstract}
We are interested in numerical schemes for the simulation of large scale gas networks. Typical models are based on 
the isentropic Euler equations with realistic gas constant. The numerical scheme is based on transformation of conservative
variables in Riemann invariants and its corresponding numerical dsicretization.  A particular, novelty of the proposed method
is the possbility to allow for an efficient discretization of the boundary and coupling conditions at nodal points of the network. 
The original discretization is analysed in view of its property to correctly recover steady states as well as to resolve
possible analytic solutions. Comparisons with existing methods show the advantage of the novel method.
\end{abstract}
%%%%%%%%%%%%%%%%%%%%%%%%%%%%%%%%%%%%%%%%%%%%%%%%%%%%%%%%%%%%%% 
\section{Introduction}

Mathematical models for transport of high--pressure  natural gas through a pipeline system has been subject of active research both in engineering literature e.g. \cite{Osiadacz1989aa,Backhaus2015aa,Andersson2016aa,Simone2019aa,PSIG2019aa} as
well as in mathematical literature, see e.g. \cite{BandaHertyKlar2006aa, BandaHertyKlar2006ab,KolbLangBales2010aa,HertyMohringSachers2010aa,BrouwerGasserHerty2011aa, ColomboGuerraHerty2009aa,Reigstad2015aa,GugatUlbrich2017aa,Egger2018aa} and the references therein. Depending on the application and physical regimes
different mathematical descriptions can be used for transport in pipe systems. This leads to a hierarchy of models available and we refer e.g. to \cite{HertyMohringSachers2010aa,RufflerMerhmannHante2018aa} for further details.
Besides the mathematical models for the gas flow in the pipe different conditions for coupling the flow at pipe--to--pipe intersections have been proposed and we refer to 
\cite{BandaHertyKlar2006aa, BandaHertyKlar2006ab,Reigstad2015aa} for modeling aspects as well as to \cite{ColomboMauri2008aa,ColomboHertySachers2008aa,ColomboGaravello2006aa}
for well--posedness results. 
\par 
In high--pressure and long--distance pipelines typical pressure mass flux suggests to neglect inertia and gravity effects in the mathematical model  \cite{HertyMohringSachers2010aa,BrouwerGasserHerty2011aa,Andersson2016aa} . Those models are also called friction dominated models and they can be obtained through asymptotic analysis. They are independent of temperature. The governing equations are given by 
\begin{equation}\label{simplified}
\partial_t \rho(t,x) + \partial_x ( (\rho u)(t,x) ) = 0, \; \partial_x p(t,x) = - \frac{f_g}{2d} (\rho u)(t,x) |u(t,x)| .
\end{equation}
Here, $\rho(t,x)$ denotes the gas density at time $t\geq 0$ and position $x \in I$ where $I=[0,L]$ and $L$ is the length of the pipe. The factor $f_g$ is the friction factor
and $d$  the diameter of the pipe.
%the cross-section of the pipe of diameter $d.$ 
The gas velocity is denoted by $u(t,x)$ and  the pressure by $p(t,x)$.% and the mass flux by $q(t,x)=a(\rho u)(t,x).$ 
The previous equation is not
closed and density--pressure relation referred to as pressure law needs to be prescribed. For isentropic Euler equations the relation is given by $p = C\rho^\gamma$ and the value of $\gamma$ for ideal gas is $\gamma=1.4.$ 
 Below we discuss further choices in detail.  If inertia effects are accounted for, equation \eqref{simplified} reads 
\begin{equation}\label{simplified2}
\partial_t \rho(t,x) + \partial_x (  (\rho u) (t,x) ) = 0, \;\partial_t (\rho u) (t,x) +  \partial_x \left( (\rho u^2)(t,x) + p(t,x) \right)= - \frac{f_g}{2d} (\rho u)(t,x) |u(t,x)| .
\end{equation}
Recently, the mathematical discussion has been extended to nonlinear hyperbolic models for gas flows using generalized pressure laws of the type  
\begin{align} \label{eq:gen press}
p = z(p) \rho. 
\end{align}
The function $z=z(p)$ is called compressibility factor.  In \cite{AlmeidaVelasquezBarbieri2014aa} different compressibility factors \eqref{eq:gen press} have been compared, both analytically and
by data obtained of measurements of a natural gas pipeline. Therein, a factor 
\begin{equation}\label{PST} z(p) = 1+ \alpha p \end{equation}
 for some $\alpha<0$ has been proposed. Classical solutions to the  isentropic Euler equations \eqref{simplified2} and the pressure law \eqref{PST} has been analysed in detail in \cite{GugatUlbrich2018aa}.  Also,  steady states of this system have been analysed in \cite{GugatUlbrich2017aa}. In this work we focus on a suitable numerical discretization of
 equation \eqref{simplified3} in the presence of general compressibility factors \eqref{eq:gen press} and a numerical formulation suitable to treat gas networks. 
 \par 
 In \cite{BalesKolbLang2009aa} it has been argued that for $|u|<<\sqrt{ \partial_\rho p(\rho) }$ a semilinear model can be derived. This is obtained by neglecting the term $\partial_x (\rho u^2)(t,x)$ but retain $\partial_t (\rho u).$ In the case of $p(\rho)$ given by the isothermal Euler equations, i.e., $p(\rho)=c^2 \rho,$ the obtained model \eqref{simplified2} is a linear wave equation. In the 
 case of a generalized pressure law \eqref{eq:gen press} the model is given by 
\begin{equation}\label{simplified3}
\partial_t \rho(t,x) + \partial_x \left(  \rho u (t,x) \right) = 0, \;\partial_t (\rho u) (t,x) +  \partial_x \left( p(t,x) \right)= - \frac{f_g}{2d} (\rho u)(t,x) |u(t,x)| \mbox{ and } p = z(p) \rho.
\end{equation}

This model is also studied widely in the literature, see for example \cite{RufflerMerhmannHante2018aa,BennerGrundelHimpeetal2018aa,DyachenkoZlotniketal2017aa}. Mostly however one quickly assumes also that $z(p)$ is constant.  Let $a$ be the cross-section of the pipe and define by $q(t,x)=a\rho(t,x)u(t,x)$ the mass flux, equation \eqref{simplified3} can be written as
\begin{equation}\label{simplified4}
\partial_t \rho(t,x) + \frac{1}{a}\partial_x  q (t,x) = 0, \;\frac{1}{a}\partial_t q (t,x) +  \partial_x p(t,x) = - \frac{f_g}{2da^2} \frac{q(t,x)|q(t,x)|}{\rho(t,x)} \mbox{ and } p = z(p) \rho.
\end{equation}

\section{Qualitative Properties of Model \eqref{simplified3} }

Prior to the numerical discretization we discuss some properties. Assuming that $a=1$, the system \eqref{simplified4} enjoys similar  
properties as the  $p-$system in Lagrangian coordinates.  Therefore, we do not repeat a discussion of its properties but only state the properties relevant for the numerical scheme later on. In the following we consider the function $\rho \to p(\rho)$ 
which is implicitly defined by $p=z(p)\rho.$   Further, we discuss properties in terms of the conservative variables $(\rho, q).$ The Jacobian of the flux function is given by 
$\begin{pmatrix} 0 & 1 \\ \partial_\rho p & 0 \end{pmatrix}$ 
with $\partial_\rho p = z(p) / ( 1 - \rho z'(p) )$. The eigenvalues are 
 \begin{equation}
\lambda^+(\rho) =  \sqrt{ \partial_\rho p(\rho)} \mbox{ and } \lambda^-(\rho) = -  \sqrt{ \partial_\rho p(\rho)}.
\end{equation}
To obtain strict hyperbolicity we impose the assumption 
\begin{equation}\label{ass1} 
\partial_\rho p(\rho) > 0 \; \forall \rho>0.
\end{equation}
 The assumption \eqref{ass1} is fulfilled  if we assume  $z'(p) \leq 0$. Note that 
 in the case of the isentropic Euler equations  and for a pressure law of the type \eqref{PST} the assumption 
 \eqref{ass1} is fulfilled.  (Right) eigenvectors to the eigenvalues $\lambda^{\pm}(\rho)$ 
 are $r^{\pm} = ( 1, \lambda^{\pm}(\rho) )^T$, respectively. 
Both characteristic fields are genuine nonlinear provided that $\partial_{\rho \rho} p \not =0.$ Note that in the case
of isothermal Euler equations ($z(p)=c^2$) the fields are linearly degenerated. As noted before   system \eqref{simplified3}  reduces to a linear wave equation. In the case of the pressure law \eqref{PST} 
 the condition  $\partial_{\rho \rho} p(\rho) \not = 0$ is fulfilled due to the assumption of strictly hyperbolicity.
 % In the general case 
 %the condition $\partial_{\rho \rho} p(\rho) \not = 0$ is fulfilled provided that equation \eqref{ass1} holds true. \textcolor{red}{I do not understand why this is true?}
 \par Since $(\lambda^+)^2=(\lambda^-)^2=\partial_{\rho}p(\rho)$ under Assumption \eqref{ass1} we have 
 \begin{equation}
 \partial_t p=(\lambda^+)^2\partial_t\rho \label{eq:rhop}
 \end{equation}
 Provided \eqref{ass1} holds true, the system is a $2\times 2$ hyperbolic balance law and therefore
 two Riemann invariants exists, in the following denoted by ${w}^\pm$ and given by 
 \begin{equation}\label{RI}
 {w}^+(t,x) = \frac{1}2 \left( q + \int_0^\rho \lambda^+(s) ds \right), \; {w}^-=
  \frac{1}2 \left( q + \int_0^\rho \lambda^-(s) ds \right). 
  \end{equation}
  The Riemann invariants are transported with speed $\lambda^\pm,$ i.e.,  ${w}^\pm$ fulfills 
  \begin{equation}\label{upwind}
  \partial_t {w}^\pm(t,x)  + \lambda^\pm(\rho(w^+,w^-)(t,x)) \partial_x {w}^\pm(t,x) = - \frac{1}2 \frac{f_g}{ 2d} (\rho u)(w^+,w^-)(t,x) |u(w^+,w^-)(t,x)|.
  \end{equation}
Here, $\rho(w^+,w^-)$ and $u(w^+,w^-)=\frac{q}{\rho}(w^+,w^-)$ are density and velocity obtained by inverting equation \eqref{RI}. The precise formulas will be given below. Under the assumption \ref{ass1} , with $\rho>0$ we have
\begin{equation}\label{ass2} 
\lambda^-(\rho) < 0 < \lambda^+(\rho),
\end{equation}
and for smooth solutions the values $w^\pm$ are transported along characteristics with slope $\lambda^\pm,$ respectively.  Also, shocks (or contact discontinuities in the case $p(\rho)=c^2\rho$) are not observed under typical gas operational conditions \cite{BalesKolbLang2009aa,HertyMohringSachers2010aa,GugatUlbrich2018aa}.  The equations \eqref{upwind} are  in non--conservative form. It is therefore not well--defined in the case of discontinuous $w^\pm.$ However, this form  exhibits the transport nature of the problem and we propose a numerical 
discretization related to the transport character of equation \eqref{upwind} for sufficiently smooth solutions.  This discretization allows to identify  correct boundary conditions  that may not be obtained in the case of applying a central discretization schemes towards the conservative formulation of the problem.  
\par
Finally, we note that equation \eqref{simplified4} allows for explicit solutions that can be used for 
validation of the numerical scheme. The solutions are found similarly to the approach in \cite{GugatUlbrich2017aa}:
Let $\rho(x,t)=\rho_0$ be any positive constant. Then, $q(t,x)=(\rho u)(t,x)=\frac{1}{C_0 + C_1 t}$ is a solution to equation \eqref{simplified4}.  Indeed, $p/z(p)=\rho_0$ is constant and implies $q$ is independent of $x.$ Therefore, conservation of mass holds true. 
Furthermore, 
\begin{equation} \partial_t q = - q^2  C_1 , \end{equation}
and therefore the conservation of momentum is satisfied for $C_1 = \frac{f_g}{2d \rho_0}$ and postive $q$. We have $C_1>0$ and if $C_0>0$, then $q\geq 0$ for $t\geq0$ and 
$(\rho_0, q)$ is a solution.  $C_0$ is a degree of freedom in the solution that can be used to match possible boundary conditions. 
\par
In \cite {GugatSchultzWintergerst2016aa} traveling wave solutions to equation \eqref{simplified2} and \eqref{PST} have been studied.  For equation \eqref{simplified4} we can proceed in a similar fashion. Consider a pressure law of the type \eqref{PST} with $\alpha \equiv -1.$ Then, we obtain an equation for the pressure $p(t,x)=g(t,x)$ as follows 
\begin{equation} 
p(t,x) := g(t,x), \; z(p)=1-p, \rho(t,x)=g(t,x)/(1-g(t,x)).
\end{equation}
A traveling wave solution $g$ is of type  $g(t,x) = y(c_1 t - c_2 x).$ For any $g$ of this type we have by definition
 $\partial_t g(t,x)  + \partial_x g(t,x) =0$ provided that $c_1=c_2.$ With $q=\rho=  \frac{g}{1-g}$ and $c_1=c_2$, we obtain 
 conservation of mass
 \begin{equation} \partial_t \rho(t,x) + \partial_x q(t,x)=0.
 \end{equation}
An equation for $g$ is obtained from the momentum equation. Written in terms of $y=y(\cdot)$ and for $C=  \frac{ f_g}{ 2d c_2}$  it reads  
\begin{equation}
y'(s)   ( 1- (1-y(s))^2) = - C y(s) ( 1- y(s)).
\end{equation} 
Its explicit solution fulfills $(y-1) \exp(y) = C \;t  + y(0)$ 
and therefore a closed form using  Lambert--W function can be given. 

%%%%%%%%%%%%%%%%%%%%%%%%%%%%%%%%%%%%%%%%%%%%%%%%%
\section{Numerical Discretization}

To derive the scheme 
we consider the system \eqref{simplified4} in a single pipe. The pipe is paramterized by $x \in [0,1]$ and
we assume \eqref{ass1} holds true.  
The spatial domain is discretized in $n$ equidistant intervals of size $\Delta x.$ The center of each cell
is denoted by $x_i = i \Delta x$, $i=0,\dots,n$ and we assume $n$ is such that 
$x_{n}= 1.$  We are interested in a discretization for $(\rho,q)$ based on a discretization 
of the formulation in Riemann invariants. For $a\not =1$ the Riemann invariants are given by 
\begin{equation}\label{RI disc}
w^{\pm}=\frac{1}{2}\left( \frac{q}{a} + \int_0^{\rho} \lambda^\pm(s)ds \right)
\end{equation}
where $\lambda^\pm(\rho)=\pm \sqrt{\frac{\partial p(\rho)}{\partial \rho}}.$ We further denote
by $f(\rho,q)=-\frac{f_g}{2da^2}\frac{q|q|z(p)}{p}$ and $p=z(p)\rho.$  Assuming sufficiently smooth 
solutions equation \eqref{upwind} will be discretized using a first--order finite--volume scheme. 
We denote the local cell average in cell $C_i$ as  $C_i=[x_{i-\frac{1}2}, x_{i+\frac12}]$ by  $w_i^\pm$  
\begin{equation}
w_i^\pm(t) = \frac{1}{\Delta x} \int_{x_{i-\frac{1}2}}^{x_{i+\frac{1}2}} w^\pm (t,\xi)d\xi. 
\end{equation}
The piecewise constant reconstruction by cell $C_i$ with cell average $w_i^\pm$ is denoted by $\tilde{w}^\pm(t,x),$ i.e., $\tilde{w}^\pm(t,x)=w_i^\pm \chi_{C_i}(x).$ 
 For simplicity of notation we still use  $(\rho,q)$   computed using Riemann invariants $(w^+,w^-)$. Integration yields 
\begin{equation}\label{upwind2a}
\partial_t w_i^\pm(t) + \frac{1}{\Delta x} \left(  \int_{x_{j-\frac{1}2}}^{x_{j+\frac{1}2}} \lambda^\pm(\rho)  
\partial_x w^\pm(t,\xi) d\xi \right) = \frac{1}{2\Delta x} \int_{x_{j-\frac{1}2}}^{x_{j+\frac{1}2}}  f(\rho,q)(t,\xi) d\xi. 
\end{equation}  
Within cell $C_i$ we use a constant reconstruction of the functions $\tilde{\rho}(t,x)$ and $\tilde{q}(t,x)$ by
the cell averages $w_i^\pm(t).$ The corresponding cell average is denoted by $\rho_i(t)$ 
for $i=0,\dots,n.$ Approximating $\rho(t,x) = \tilde{\rho}(t,x)$  and 
 a midpoint rule to the source term 
we obtain up to order $O(\Delta x^2)$ 
\begin{equation}\label{upwind2}
\partial_t w_i^\pm(t) + \frac{ \lambda^\pm(\rho_i( t) )  }{\Delta x} 
\left( \tilde{w}^\pm( t, x_{i+\frac{\Delta x}2} ) - \tilde{w}^\pm( t, x_{i-\frac{\Delta x}2} )
   \right) =  \frac12 f(\rho_i(t),q_i(t)).
\end{equation}  
Due to assumption \eqref{ass2} we use an Upwind discretization. This implies that  
since $\lambda^+>0$ we approximate  $\tilde{w}^+( t, x_{i+\frac{\Delta x}2} )=w^+_i$
for $i=0,\dots,n,$ and  similarly since $\lambda^-<0$ we approximate 
 $\tilde{w}^-( t, x_{i+\frac{\Delta x}2} )=w^-_{i+1}.$ This leads to the 
 following discretization for $f_i(t)=f(\rho_i(t),q_i(t))$ 
\begin{align}
\dt w_i^+(t)+ \frac{\lambda^+(\rho_i(t))}{\Delta x}
 \left( w^+_{i}(t) - w^+_{i-1}(t) \right)  &=\frac12 f_i(t) \label{eq:wminus}\quad i=1,\dots n\\
\dt w_i^-(t)  +\frac{\lambda^-(\rho_i(t))}{\Delta x}   \left(  w^-_{i+1}(t)-w^-_{i}(t) \right)
 &=\frac12  f_i(t) \label{eq:wplus}\quad i=0,\dots n-1
\end{align}
Any explicit  temporal discretization  on $(t_m, t_m+\Delta t), m>0$ with time step $\Delta $ 
  has to fulfill the CFL condition 
\begin{equation}
\Delta t \leq \frac{ \Delta x}{ \max\limits_{i=0,\dots,n}  \lambda^+(\rho_i(t_m))  }.
\end{equation}
Equations \eqref{eq:wminus} and \eqref{eq:wplus} require boundary conditions for $w^+_0(t)$ at $x=x_0$ 
and boundary conditions for $w^-_n(t)$ at $x=x_n,$ respectively.
\medskip
\paragraph{Discretization in conservative variables $(\rho,q)$} 
The numerical scheme  \eqref{eq:wminus} and \eqref{eq:wplus} is reformulated in terms of cell averages 
of conservative variables $(\rho_i,q_i)(t).$ Further, we denote by 
\begin{equation} \lambda_i = \lambda^+(\rho_i) \end{equation}
and we have $\lambda^-(\rho_i)=-\lambda_i.$ We compute 
\begin{align}
\lambda_i \dt\rho_i=&\dt(w_i^+-w_i^-)=-\frac{\lambda_i}{\Delta x}(w_{i}^+-w_{i-1}^+-w_i^-+w_{i+1}^-)+\frac12 \left( f_i -f_i \right)\\
%=& -\frac{\lambda_i}{2 \Delta x}\left(\frac{q_{i}}{a}+\int_0^{\rho_{i}} \lambda(s)ds-\frac{q_{i-1}}{a}-\int_0^{\rho_{i-1}} \lambda(s)ds-\frac{q_i}{a}+\int_0^{\rho_i} \lambda(s)ds+\frac{q_{i+1}}{a}-\int_{0}^{\rho_{i-1}} \lambda(s)ds \right)\\
=& -\frac{\lambda_i}{2 \Delta x}(\frac{q_{i+1}-q_{i-1}}{a}-\int_{\rho_{i}}^{\rho_{i+1}} \lambda(s)ds + \int_{\rho_{i-1}}^{\rho_{i}} \lambda(s)ds)\\
=& -\frac{\lambda_i}{2 \Delta x}\left(\frac{q_{i+1}-q_{i-1}}{a}-\frac{\rho_{i+1}-\rho_{i}}{2}( \lambda_{i+1}+\lambda_{i})+\frac{\rho_{i}-\rho_{i-1}}{2}( \lambda_{i}+\lambda_{i-1})\right ) +O(\Delta x)\\
=&-\frac{\lambda_i}{2 \Delta x}\left(\frac{q_{i+1}-q_{i-1}}{a} \right ) +O(\Delta x),
\end{align}
and similarly 
\begin{align}
\frac{1}{a}\dt q_i =&\dt (w_i^+ + w_i^-)=\frac{\lambda_i}{\Delta x}(w_{i-1}^+-w_i^++w_{i+1}^--w_{i}^-)+ \frac12 \left( f_i+f_i \right) \\
%=& -\frac{\lambda_i}{2 \delta x}(q_{i+1}+\int_0^{\rho_{i+1}} \lambda(s)ds-q_i-\int_0^{\rho_i} \lambda(s)ds+q_i-\int_0^{\rho_i} \lambda(s)ds-q_{i-1}+\int_{0}^{\rho_{i-1}} \lambda(s)ds)+f_i\\
=& \frac{\lambda_i}{2 \Delta x}\left(\frac{q_{i+1}-2q_i+q_{i-1}}{a}-\int_{\rho_{i-1}}^{\rho_{i+1}} \lambda(s)ds\right)+f_i\\
%=& -\frac{\lambda_i}{2 \delta x}\left (q_{i+1}-q_{i-1}+\frac{\rho_{i+1}-\rho_i}{2}( \lambda_i+\lambda_{i+1})-\frac{\rho_{i}-\rho_{i-1}}{2}( \lambda_i+\lambda_{i-1})\right)+f_i\\
=&- \frac{\lambda_i}{2 \Delta x}\left (\frac{\rho_{i+1}-\rho_{i-1}}{2}( \lambda_{i+1}+\lambda_{i-1}))\right)+f_i + 
O(\Delta x)
\end{align}

Hence, a semi--discretization for the cell averages in conservative variables is given by 
\begin{align}
\dt p_i(t)=&-\frac{\lambda_i^2}{2 \Delta x a}\left(q_{i+1}(t)-q_{i-1}(t)\right )\quad i=1,\dots n-1\\
\dt q_i(t) =& -\frac{\lambda_i a}{4 \Delta x}\left (\frac{p_{i+1}(t)}{z_{i+1}(t)}-\frac{p_{i-1}(t)}{z_{i-1}(t)}\right)( \lambda_i+\lambda_{i+1})+a f_i \quad i=1,\dots n-1
\end{align} 
where $p_i(t)=p(\rho_i(t))$ and $z_i(t)=z(p(\rho_i(t)))$ and using \eqref{eq:rhop} and $p=z(p)\rho$. The previous equation has to be complemented with 
suitable boundary conditions obtained by equation \eqref{eq:wminus} for $i=0$ and by equation 
\eqref{eq:wplus} for $i=n,$ i.e., 
\begin{align}\label{eq:w bc}
\dt w_n^+(t)= -\lambda_n\frac{w^+_{n}(t)-w^+_{n-1}(t)}{\Delta x}+ \frac12 f_n \mbox{ and } 
\dt w_0^-(t)= \lambda_0\frac{w^-_{1}(t)-w^-_{0}(t)}{\Delta x}+ \frac12 f_0.
\end{align}
Using equation \eqref{eq:rhop}
%\begin{equation} 
%\dt p=\frac{\partial p}{\partial \rho}\dt \rho=\lambda(p)^2\dt \rho, \: 
%\lambda\dt\rho=\frac{1}{\lambda}\dt p \end{equation}
those equations allow to derive 
boundary conditions for $(\rho,q)_i$, $i\in \{ 0, n \}$:
\begin{align}
\partial_t \frac{q_n}{a}+\frac{1}{\lambda_n}\partial_t p_n &=-\lambda_n\frac{q_n-q_{n-1}}{a \delta x}-\frac{\lambda_n}{2\delta x}(\frac{p_n}{z_n}-\frac{p_{n-1}}{z_{n-1}})(\lambda_n+\lambda_{n-1})+f_n\label{bndeq1}\\
\partial_t \frac{q_0}{a}-\frac{1}{\lambda_0}\partial_t p_0 &=\lambda_0\frac{q_1-q_{0}}{a\delta x}-\frac{\lambda_0}{2\delta x}(\frac{p_1}{z_1}-\frac{p_0}{z_0})(\lambda_1+\lambda_0)  +f_0\label{bndeq2}
\end{align}
Additionaly, we assume initial conditions $p_{IC}(x), q_{IC}(x)$ and boundary conditions $p_{BC}(t), q_{BC}(t)$ given: 
\begin{equation}
p(0,x)=p_{IC}(x), \; q(0,x)=q_{IC}(x), \; p(0,t)=p_{BC}(t) \mbox{ and } q(1,t)=q_{BC}(t).
\end{equation}
For a single pipe we therefore obtain the following semi--discretized system together with discretized 
initial and boundary conditions 
\begin{align}
\dt p_i=&-\frac{\lambda_i^2}{2 \Delta x a}\left(q_{i+1}-q_{i-1}\right ), \; p_i(0)=p_{IC,i}, \; i=1,\dots,n-1\\
\dt q_i =& -\frac{\lambda_i a}{4 \Delta x}\left (\frac{p_{i+1}}{z_{i+1}}-\frac{p_{i-1}}{z_{i-1}}\right)( \lambda_i+\lambda_{i+1})+af_i, \; q_i(0)=q_{IC,i},  \; i=1,\dots,n-1 \\
\partial_t \frac{q_n}{a}+\frac{1}{\lambda_n}\partial_t p_n=&-\lambda_n\frac{q_n-q_{n-1}}{a \Delta x}-\frac{\lambda_n}{2\Delta x}(\frac{p_n}{z_n}-\frac{p_{n-1}}{z_{n-1}})(\lambda_n+\lambda_{n-1})+f_n, \; p_n(0)=p_{IC,n}\\
\partial_t \frac{q_0}{a}-\frac{1}{\lambda_0}\partial_t p_0=&\lambda_0\frac{q_1-q_{0}}{a\Delta x}-\frac{\lambda_0}{2\Delta x}(\frac{p_1}{z_1}-\frac{p_0}{z_0})(\lambda_1+\lambda_0)  +f_0, \; q_0(0)=q_{IC,0} \\
p_0(t) &= p_{BC}(t), \; q_n(t)=q_{BC}(t) .  \label{bc01}
\end{align}
%%%%%%%%%%%%%%%%%%%%%%%%%%%%%%%%%%%%%%%%%%%%%%%%%%%%%%5
\paragraph{Other discretization schemes}
In the numerical experiments we compare this discretization against others from the literature. The first one is given in  \cite{GruJHetal14} and uses a midpoint discretization:

\begin{align}
\frac{1}{2\lambda_{i+1}^2}\dt p_{i+1}+\frac{1}{2\lambda_{i}^2}\dt p_{i}=&-\frac{1}{\Delta x a}\left(q_{i+1}-q_{i}\right ),\; i=0,\dots,n-1\\
\dt \frac{q_i+q_{i+1}}{2} =& -\frac{ a}{\Delta x}\left ( p_{i+1}-p_i\right)-\frac{f_g}{4da}\frac{(q_i+q_{i+1})|q_i+q_{i+1}|}{p_i+p_{i+1}},   \; 0=1,\dots,n-1 \\
q_i(0)=q_{IC,i},  \; 0=1,\dots,n-1& \;p_i(0)=p_{IC,i}, \; i=1,\dots,n\;p_0(t) = p_{BC}(t), \; q_n(t)= q_{BC}(t). 
\end{align}
Another one can be found in \cite{GruJ15,StoM18} and will be called endpoint discretization in the following:

\begin{align}
\dt p_i=&-\frac{\lambda_i^2}{ \Delta x a}\left(q_{i}-q_{i-1}\right ), \; p_i(0)=p_{IC,i}, \; i=1,\dots,n\\
\dt q_i=& -\frac{ a}{\Delta x}\left ( p_{i+1}-p_i\right)-\frac{f_g}{2da}\frac{q_i|q_i|}{p_{i+1}}, \; q_i(0)=q_{IC,i},  \; i=0,\dots,n-1 \\
p_0(t) &= p_{BC}(t), \; q_n(t)= q_{BC}(t). 
\end{align}

\medskip
\par
\paragraph{Consistent discretization of steady--states} 

Schemes that preserve steady states exactly are called well-balanced, and their development is a lively topic in the field of hyperbolic balance laws, see e.g. \cite{GreenbergLeroux1996aa,AudusseBouchutBristeau2004aa,BollermannChenKurganov2013aa,NoellePankratzPuppo2006aa,ChertockCuiKurganov2015aa,XingShu2013aa} and references therein. Usually, these schemes use specific knowledge of  an equilibrium state. For the proposed scheme 
and in the case $z(p)=c^2$ we obtain the following results:  The proposed scheme conserves the continuous steady state at most to order $O(\Delta x)$ and the scheme conserves discrete steady states exactly. 
\par 
Provided $z(p)=c^2$ for some constant $c.$  Then, the continuous steady state of the system \eqref{simplified3} are 
\begin{equation}\label{css} q(x)= C_q \mbox{ and }  \partial_x p(x)=  f( c^2 p(x), C_q) \end{equation}
where we recall  $ f(\rho,q)=- \frac{f_g}{2d a^2  } \frac{ q |q| }{ \rho }.$ If we assume that 
the data $p=p(x)$ is given by equation \eqref{css} then, the discretized steady states are at the cell center $x_i$ are given by 
\begin{equation} \label{eq:c ss} 
p(x_i)=\int_0^{x_i} f( c^2 p(y),C_q) dy \mbox{ and } q(x_i)= C_q.
\end{equation} Applying the previous scheme with initial conditions 
given by equation \eqref{eq:c ss}  yields 
 \begin{align} 
 \partial_t p_i(t) &= 0, \\
\partial_t q_i(t) &= - \frac{ c a }{4 \Delta x} \left( \frac{1}{c^2} \int_{x_{i-1}}^{x_{i+1}} f(c^2 p(y),C_q) dy \right) ( 2c ) + a f_i  = - \frac{a}{2\Delta x} \int_{x_{i-1}}^{x_{i+1}} f(c^2 p(y),C_q) dy + a f_i.
 \end{align}
By definition  we have $f_i=f(\rho_i,q_i)$ and therefore the last equation is equal to zero up to the order $\Delta x.$ Hence, 
$\partial_t q_i(t) =  O(\Delta x)$ and the scheme approximates the continuous steady state up to order $O(\Delta x).$ 
Clearly, it is possible to define higher--order integration of the source term, e.g.,  replacing $f_i$ by $\sum \omega_j f_{j}$ 
where $\omega_j$ are integration weights and $f_j=f(q_j,\rho_j)$ with $j \in \{ i-1,i,i+1\}.$ By the proposed Upwind scheme 
it is not reasonable to have stencils beyond $i-1$ and $i+1$ for the integration of the source term. This leaves
as integration schemes Newton-Cotes formulas. With three points we may use Simpson's formula with an error
of $(\Delta x)^5.$ This yields consistency in this case up to order $O(\Delta x)^4$ if we discretize
\begin{equation} f(\rho,q)=  \left( \frac{1}6 f_{i-1} + \frac{2}3 f_i + \frac{1}6 f_{i+1} \right).
\end{equation}

\par
We also consider the conservation of discrete steady statesfor the same choice of $z(p)=c^2.$ Then, $\lambda^\pm(\rho)=\pm c$ and the 
discretization simplifies to  
\begin{equation}\label{eq:d ss}
 \partial_t p_i(t) = - \frac{c^2}{2 \Delta x a } (q_{i+1}(t)-q_{i-1}(t), \; \partial_t q_i (t)= - \frac{a}{2\Delta x} ( p_{i+1}(t)-p_{i-1} (t)) -  \frac{f_g c^2 q_i(t) |q_i|(t) } { 2d a p_i(t) }.
\end{equation}
Consider an explicit Euler discretization as temporal discretization and denote by $q_{k,i}=q_i(t_k), p_{k,i}=p_i(t_k).$ 
 Assume that $q_{k,i}=C_q$ for all $i$ at time $t_k$ and  $p_{k,i+1}=p_{k,i-1} - \frac{  \Delta x \; f_g c^2 C_q |C_q| } { d  p_{k,i} }.$ 
 Those are the discrete steady states of the scheme \eqref{eq:d ss}. Then, an explicit Euler discretization yields
 \begin{equation}
 p_{k+1,i} = p_{k,i} + 0, \; q_{k+1,i} = q_{k,i}  + 0 = C_q.
 \end{equation}
Therefore, the scheme preserves the discrete steady states exactly for all future times $k.$ 
%%%%%%%%%%%%%%%%%%%%%%%%%%%%%%%%%%%%%%%%%%%%%%%%%%%%%%5
\medskip
\par
\paragraph{Discretization of pipe networks}

A pipe networks is modelled as directed graph $(E,V)$ where $E$ is the set of all edges $k\in \{1,\dots,K\}$ and $V$ the set of all vertices. 
For the vertices we distinguish between internal nodes in the network $V_0$ (having a degree larger or equal to two) and boundary nodes of the network, i.e., vertices with  degree one. For vertices of degree equal to one we may either prescribe pressure or mass flux conditions 
as in equation \eqref{bc01}. For the internal nodes we use coupling conditions given by equal pressure and the conservation of 
mass \eqref{cond}. For further details and a discussion of the  coupling conditions we refer to \cite{Osiadacz1989aa,BandaHertyKlar2006aa} and references therein.
\par
We introduce the following notation.  For simplicity we assume all pipes are parameterized by $x \in [0,1]$. By $p_i^k(t), q_i^k(t)$ we denote the cell average of the pressure and mass flux in pipe $k$ and cell $C_i$ of pipe $k$ at time $t.$ Hence, the equations for the temporal 
evolution of the cell averages $(p_i^k,q_i^k), i=0,\dots,n$ and $k\in E$ are then given by equations 
\begin{align}
\dt p_i^k=&-\frac{\lambda_i^2}{2 \Delta x a}\left(q_{i+1}^k-q_{i-1}^k\right )\quad i=1,\dots n-1\label{eq1}\\
\dt q_i^k =& -\frac{\lambda_i a}{4 \Delta x}\left (\frac{p_{i+1}^k}{z_{i+1}^k}-\frac{p_{i-1}^k}{z_{i-1}^k}\right)( \lambda_i+\lambda_{i+1})+af_i^k \quad i=1,\dots n-1\\
\partial_t \frac{q_n^k}{a}+\frac{1}{\lambda_n}\partial_t p^k_n=&-\lambda_n\frac{q^k_n-q^k_{n-1}}{a \Delta x}-\frac{\lambda_n}{2\Delta x}(\frac{p^k_n}{z^k_n}-\frac{p^k_{n-1}}{z^k_{n-1}})(\lambda_n+\lambda_{n-1})+f^k_n\\
\partial_t \frac{q^k_0}{a}-\frac{1}{\lambda_0}\partial_t p^k_0=&\lambda_0\frac{q^k_1-q^k_{0}}{a\Delta x}-\frac{\lambda_0}{2\Delta x}(\frac{p^k_1}{z_1}-\frac{p^k_0}{z^k_0})(\lambda_1+\lambda_0)  +f^k_0\label{eq2}
\end{align}
and  initial conditions.  For a vertex $v\in V$ denote the set of all edges $k \in E$ incoming to $v$ by $\delta^-_v$ and all 
edges $\ell \in E$ exiting from $v$ by $\delta^+_v.$ In case $ |\delta^+_v| + |\delta^-_v| \geq 2$  the coupling conditions at the vertex $v\in V$ 
\begin{align}\label{cond}
 p_n^k(t) &=p_0^\ell(t) , \; \forall k \in \delta^-_v, \; \forall \ell \in \delta^+_v \mbox{ and } 
 \sum\limits_{k \in \delta^-_v}  q^k_n(t)  =\sum\limits_{k\in \delta^+_v} q_0^k(t)  .
 \end{align}

%%%%%%%%%%%%%%%%%%%%%%%%%%%%%%%%%%%%%%%%%%%%%%%%%%%%%%5
Conditions \eqref{cond} together with intial conditions, boundary conditions for pressure and mass flux at nodes of degree one 
and the discretization on the network form the fully discrete scheme. Let $p(t)$ be a vector stacking all discretization points of all pipes of the pressure in one vector and the same for all mass fluxes in $q(t)$. The only points not included as they are already known are the start and end nodes which have prescribed boundary conditions. Then equation \eqref{eq1}-\eqref{eq2} can be written as

\begin{align} \underbrace{\begin{bmatrix} M_1 & 0 \\ 0 & M_2 \\ B_{11}(p(t)) & B_{12}\\B_{21}(p(t))& B_{22}\end{bmatrix}}_{E_1}\partial_t\begin{bmatrix} p(t)\\q(t)\end{bmatrix}=F(p(t),q(t),p_{BC}(t),q_{BC}(t)) \end{align}

where $p_{BC}$ and $q_{BC}$ are  the pressure and mass flux boundary conditions. Adding equation \eqref{cond} into that the full equation reads
 
\begin{align} \begin{bmatrix} E_1\\ 0\end{bmatrix}\partial_t\begin{bmatrix} p(t)\\q(t)\end{bmatrix}=\begin{bmatrix}F(p(t),q(t),p_{BC}(t),q_{BC}(t)) \\A(p(t),q(t))\end{bmatrix}\label{DAE}\end{align}

Where $A(p(t),q(t)$ collects all the algebraic equations defined in \eqref{cond}. This equation \eqref{DAE} is a nonlinear differential algebraic equation.

%%%%%%%%%%%%%%%%%%%%%%%%%%%%%%%%%%%%%%%%%%%%%%%%%%%%%%5

\section{Computational Results}

In the following we consider three networks: a pipe, the diamond benchmark also described in \cite{BennerGrundelHimpeetal2018aa} and a realistic  network mentioned in \cite{Farzaneh2016aa}. The networks and the numerical simulations  are described in the following sections.
\subsection{Simulation of Gas Flow in a Single Pipe}
\begin{figure}[htbp]
	\centering
\input{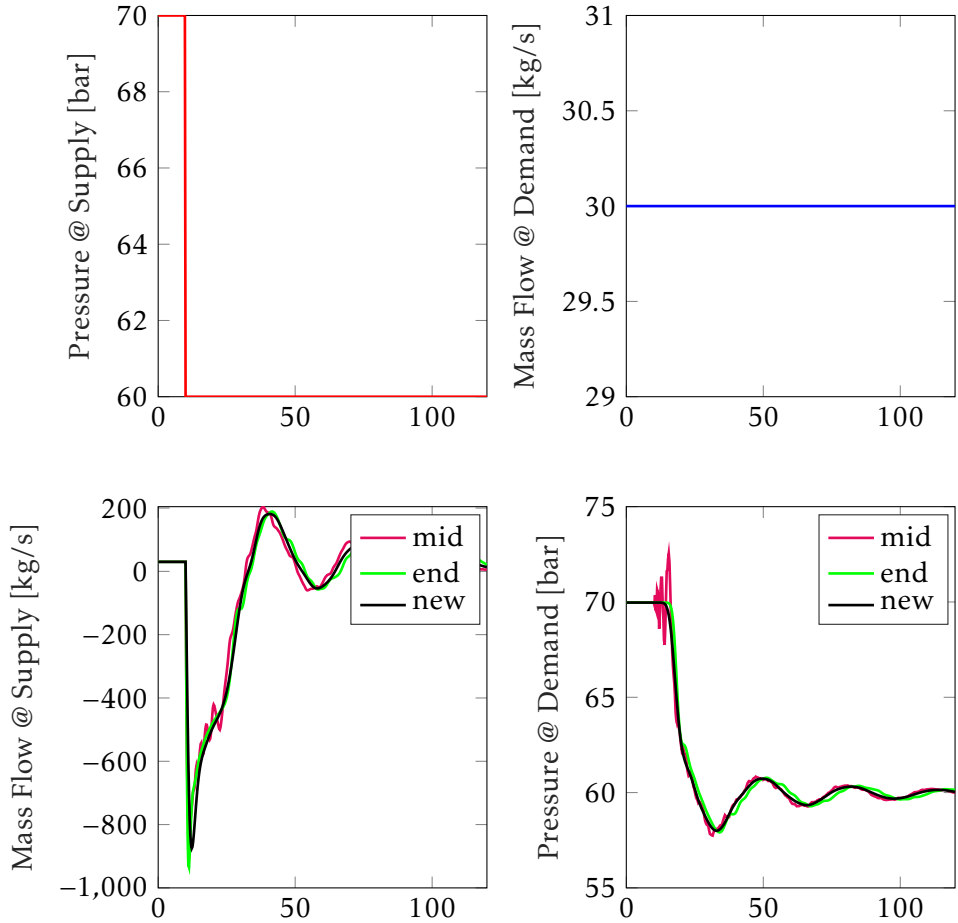}
	\caption{Numerical simulation of a pressure drop at the inlet of a pipe}
	\label{fig:pipestep}
\end{figure}

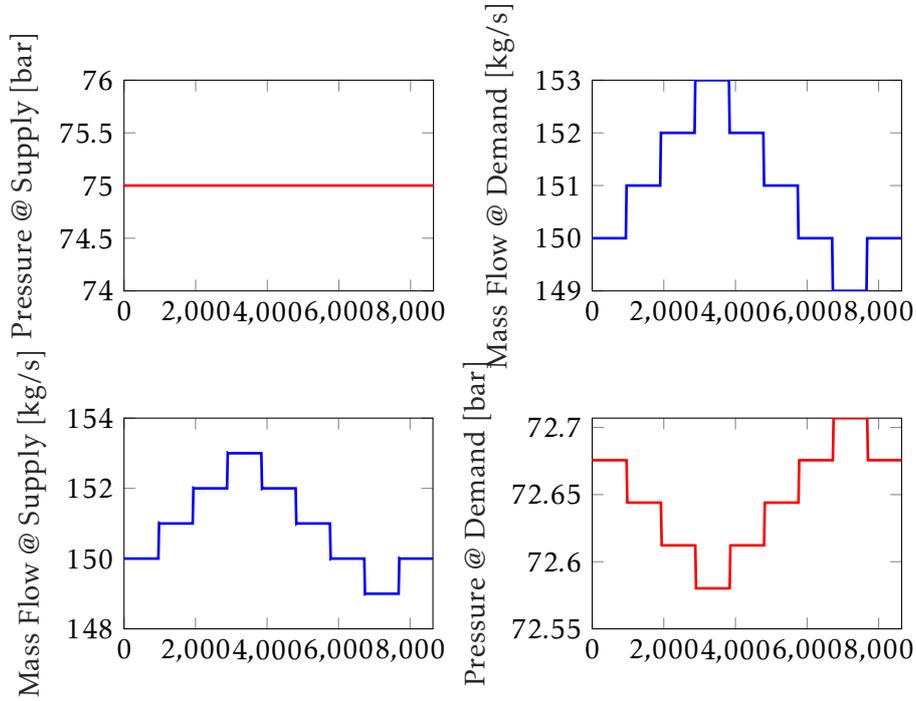
\begin{figure}[htbp]
	\centering
\input{hyp_pipe_wave}
	\caption{Wave Scenario on a pipe}
	\label{fig:pipewave}
\end{figure}

The pipe is 3 km  long and has a diameter of 0.762m and a 
%rougness of 0.0005. 
friction factor  $f_g=0.0178$ and $\lambda=383.0735$. We compare the new discretization with the two other discretizations mentioned before, one which we call midpoint discretization and abbreviated by ``mid'' and an endpoint discretization abbreviated by ``end''. In particular, consider two scenarios shown in Figure \ref{fig:pipestep} and \ref{fig:pipewave}. In the first scenario we start the system in a stationary solution and then abruptly decrease the pressure at the inlet of the pipe. We compute the pressure at the outlet and the flux at the inlet which then dynamically changes until the system goes back to a stationary solution. The flux at the outlet is kept constant. The three different results for  different discretization schemes are plotted in Figure \ref{fig:pipestep}, and we observe that the novel discretization does not lead to any oscillations in flux and pressure. The discretization `mid' produces {\bf unphysical} oscillations. 

In the second scenario we change the flux  every 1000 seconds according to Figure \ref{fig:pipewave} by keeping a constant pressure at the inlet of 75 bar. This yields  a change of the mass flow at the inlet and a change of the pressure at the outlet. Since the dynamic behaviour is not recognizable over this long time period all numerical schemes give basically the same result. 

In Table \ref{tab:time_s} we summarize the simulation time for the two scenarios named ``step'' and ``wave'' and the three different methods. Here we see that the novel method outperforms the other methods. The time is wall time and we use \texttt{ode15s} a MATLAB ODE solver, which can handle index--1 DAEs to do the time integration. 

\begin{table}[h]
    \centering
    \caption{\footnotesize Computational time (seconds) for Simulation}
    \label{tab:time_s}
    \vspace{0.3cm}
    \begin{tabular}{ccc}
        \toprule
        method  & step& wave   \\
        \midrule
        end  & 0.45& $1.27$           \\ 
        \addlinespace
        mid  & 3.98 & $202.57$          \\ 
        \addlinespace
        new  & 0.28 & $0.50$         \\ 
        \bottomrule
    \end{tabular}
\end{table}

\subsection{Resolution of Steady States} 
\begin{figure}[htbp]
	\centering
\includegraphics[width=4.5in]{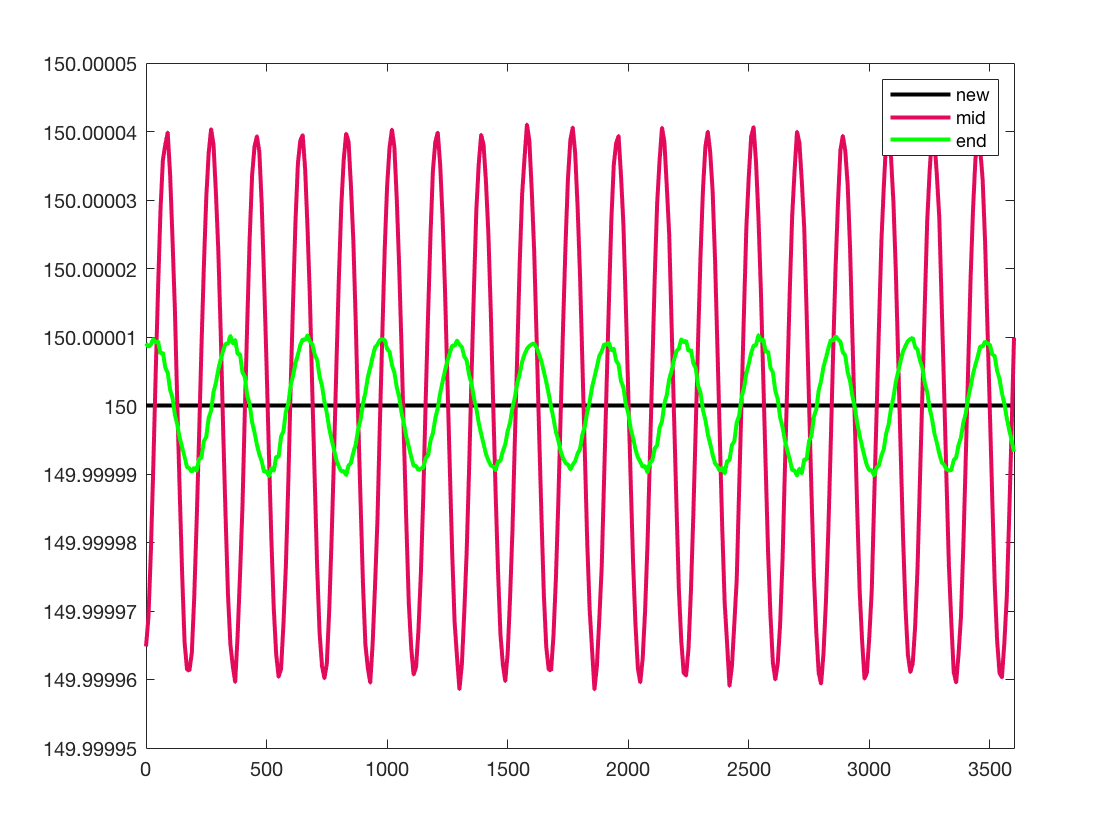}

	\caption{Oszillations of mass flux at the inlet in steady state}
	\label{fig:stat}
\end{figure}

If the boundary conditions are kept contant in time the system will converge to  steady state. This is true for all three discretizations. However the steady state of the new discretization is the only state which is actually steady. While the other  numerical states are states that yields   fluctuations around the constant steady state. For our next numerical experiment we again simulated on a pipe of length 3km, with a  diameter of 0.762m and a friction factor  $f_g=0.0178$ and $\lambda=383.0735$. We put a pressure of 155 constant on the inlet of the pipe and a constant flux of 150 kg/s on the outlet of the pipe. The stationary solution is the given by a constant mass flux in the entire pipe of 150 kg/s and a pressure of $\approx153.8887$ bar. However for the endpoint and midpoint discretization we observe oscillations that do not decay over time as shown in Figure \ref{fig:stat}.

%\begin{figure}[htbp]
%	\centering
%\input{hyp_pipe_hyp}
%
%	\caption{New discretization}
%	\label{fig:qset_cond1}
%\end{figure}
\subsection{Simulation of a Diamond--shaped Network}

The network  consists of 9 pipes, all of length 1km, diameter of 1m and a friction factor $f_g=0.0196$.  All nodes are at the same geodesic height and all pipes are flat.  We use it to show how different the three discretization schemes resolve the dynamic behaviour. 
\begin{figure}[htbp]
\centering
\includegraphics[width=4in]{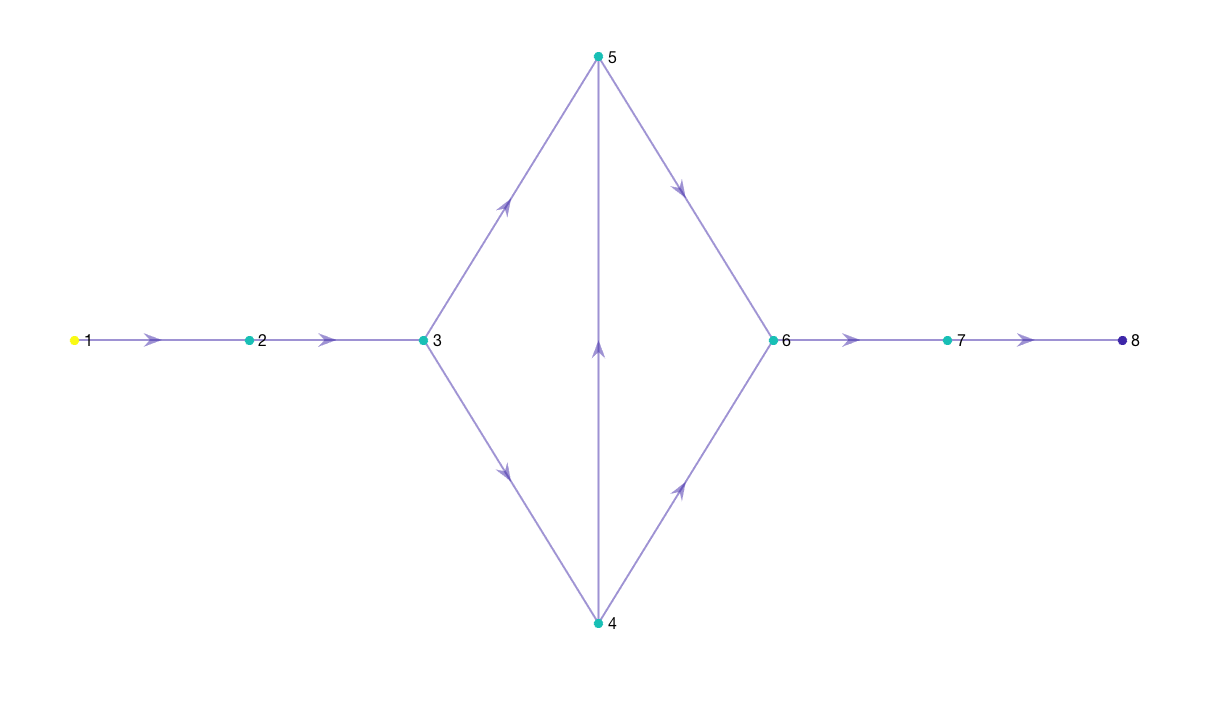}
	\caption{Topology of the diamond network}
	\label{fig:diamond_top}
\end{figure}

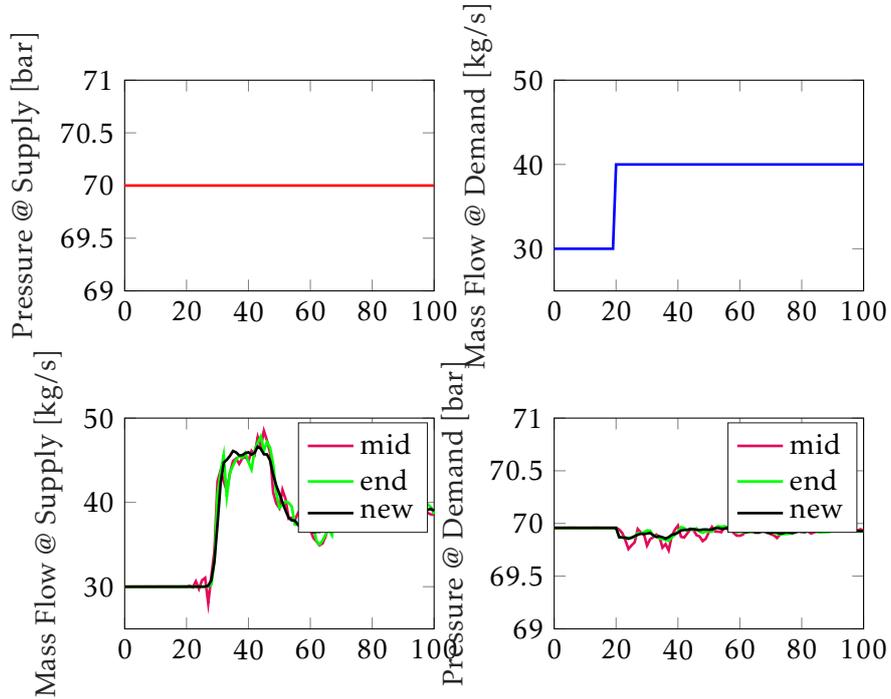
\begin{figure}[htbp]
	\centering
\input{hyp_diamond_step}

	\caption{Numerical Simulation on the diamond network}
	\label{fig:diamond_step}
\end{figure}

In this scenario we keep the pressure at the inlet constant 70 bar and increase the mass flux at the outlet of the network, the demand node from 30 kg/s to 40 kg/s. We use a constant $\lambda=673.7021$. In Figure \ref{fig:diamond_step} we see the very different restults in particular for the pressure at the supply node. Again as already seen for the simple pipe example the novel discretization does {\bf not} generate unphysical oscillations. 
%In Figure \ref{fig:diamond_updown}, we plot the pressure at node 4 and 5 of the diamond network and also the 3 fluxes at node 4 and 5. If a flux is positive it goes in the direction given in the gaph in Figure \ref{fig:diamond_top} if it is negative in the opposite direction. 

%We see that there is very little action between the two nodes and we see that on that edge direction of the flux changes. 
%\begin{figure}[htbp]
%	\centering
%\input{hyp_diamond_updown}

%	\caption{Diamond top and bottom node}
%	\label{fig:diamond_updown}
%\end{figure}

\subsection{Simulation of a Pipe Network}
\begin{figure}
	\caption{Topology of the Network}
\includegraphics[width=3in]{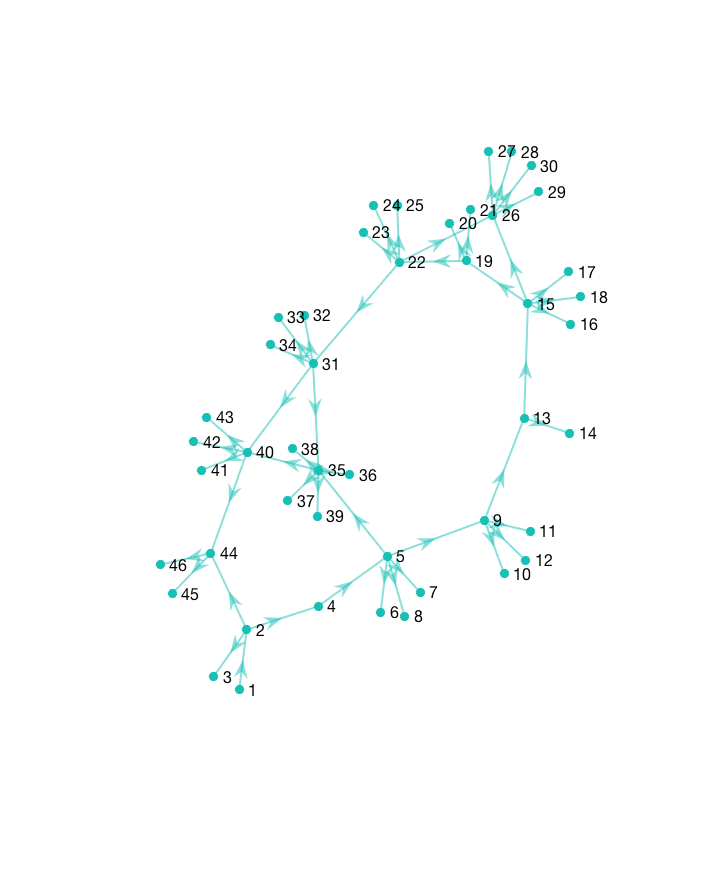}
	\label{fig:GorR15_top}
\end{figure}

\begin{figure}[htbp]
	\centering
\input{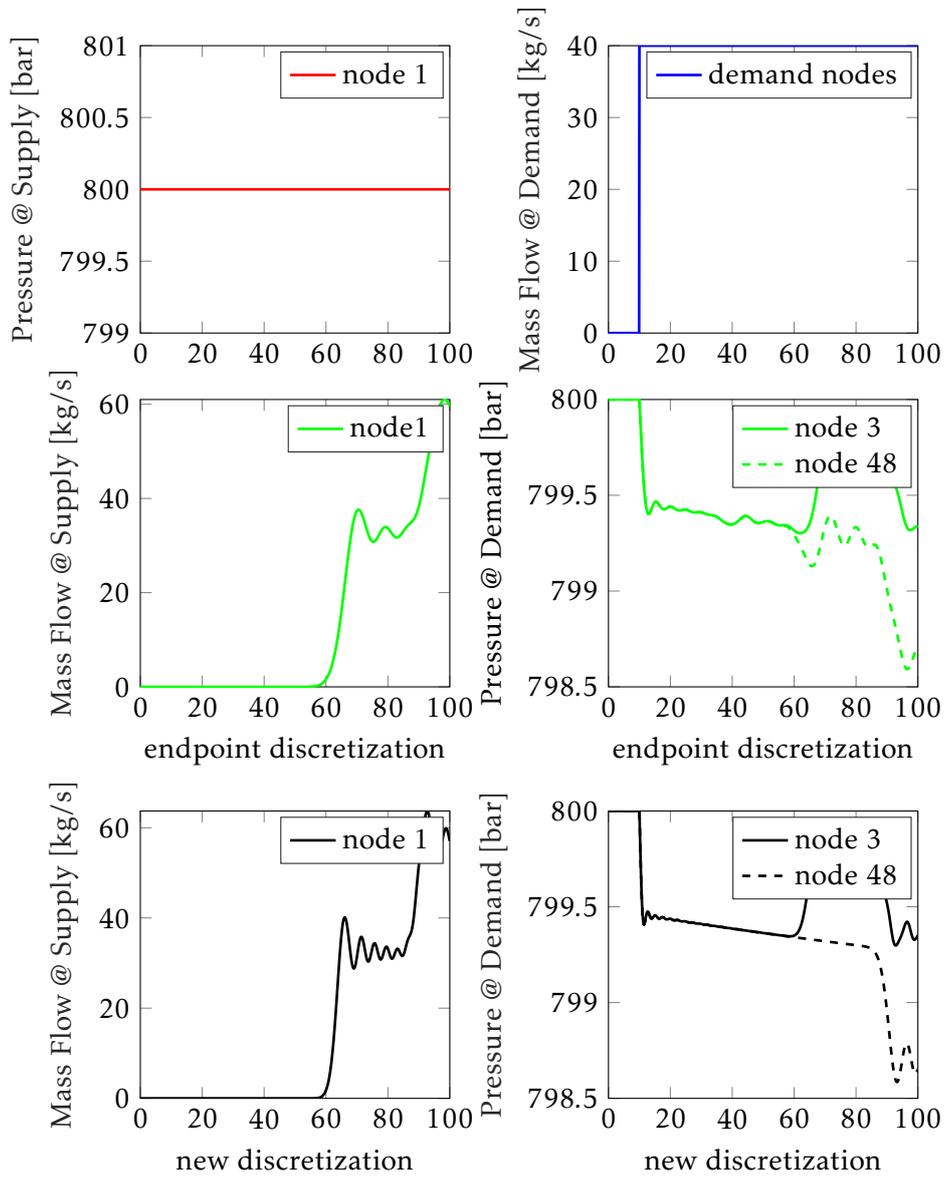}

	\caption{Simulation Result for a Network Scenario}
	\label{fig:GorR15_scn}
\end{figure}

The new discretization is also used on a   realistic network, whose topolgy is illustrated in Figure \ref{fig:GorR15_top}. Node 1 is a supply node and all the other nodes that are end nodes are demand nodes. This network has a total of  46 nodes with 1 supply nodes, 13 interior noes and 23 supply nodes. To keep it simple and reproducible the pipes are all 10 km long, 0.6 m wide and have a roughness of 0.01 resulting in a friction factor of 0.0454. Again $\lambda=383.0545$ is constant.  We simulated a scenario where  at the simple stationary point  all  fluxes are zero and all the pressures are identically 800 bar. This state is modified by  a sudden increase at the outlet flux to 40 kg/s at all demand nodes. The pressure at two different outlets and the flux at the inlet is shown in Figure \ref{fig:GorR15_scn} for the new discretization and the endpoint discretization.

\section{Conclusions}

We showed that the discretization using Riemann invariants  allows to correctly approximate the transport phenomena underlying the nonlinear model. 
The novel discretization does not lead to  any oszillations when converging to  steady state. It is also approximates the 
discrete steady state to any order. Numerical experiments show also for complex network geometries a good performance of the 
scheme. It also keeps all the advantages of the endpoint discretizations for the numerical treatment as a differential algebraic system and can be computed  efficiently.

%%%%%%%%%%%%%%%%%%%%%%%%%%%%%%%%%%%%%%%%%%%%%%%%%%%%%%%%%%%%%%%%%%%%%%%%

%%%%%%%%%%%%%%%%%%%%%%%%%%%%%%%%%%%%%%%%%%%%%%%%%%%%%%%%%%%%%%%%%%%%%%%%
\subsection*{Acknowledgment}
This work has been supported by HE5386/14,15-1, BMBF ENets 05M18PAA , BMWi mathenergy 0324019B, ERDF/EFRE: ZS/2016/04/78156 and ID390621612 Cluster of Excellence Internet of Production (IoP).
%%%%%%%%%%%%%%%%%%%%%%%%%%%%%%%%%%%%%%%%%%%%%%%%%%%%%%%%%%%%%%%%%%%%%%%%
\bibliographystyle{siam}
\bibliography{arXiv_preprint}
\end{document}

%% file: hyp_pipe_wave.tex
\begin{tikzpicture}

\begin{axis}[%
width=1.603in,
height=1.1in,
at={(1.011in,2.406in)},
scale only axis,
xmin=0,
xmax=8640,
ymin=74,
ymax=76,
ylabel style={font=\color{white!15!black}},
ylabel={Pressure @ Supply [bar]},
axis background/.style={fill=white},
legend style={legend cell align=left, align=left, draw=white!15!black}
]
\addplot [color=red, line width=1.0pt]
  table[row sep=crcr]{%
0	75\\
20	75\\
40	75\\
60	75\\
80	75\\
100	75\\
120	75\\
140	75\\
160	75\\
180	75\\
200	75\\
220	75\\
240	75\\
260	75\\
280	75\\
300	75\\
320	75\\
340	75\\
360	75\\
380	75\\
400	75\\
420	75\\
440	75\\
460	75\\
480	75\\
500	75\\
520	75\\
540	75\\
560	75\\
580	75\\
600	75\\
620	75\\
640	75\\
660	75\\
680	75\\
700	75\\
720	75\\
740	75\\
760	75\\
780	75\\
800	75\\
820	75\\
840	75\\
860	75\\
880	75\\
900	75\\
920	75\\
940	75\\
960	75\\
980	75\\
1000	75\\
1020	75\\
1040	75\\
1060	75\\
1080	75\\
1100	75\\
1120	75\\
1140	75\\
1160	75\\
1180	75\\
1200	75\\
1220	75\\
1240	75\\
1260	75\\
1280	75\\
1300	75\\
1320	75\\
1340	75\\
1360	75\\
1380	75\\
1400	75\\
1420	75\\
1440	75\\
1460	75\\
1480	75\\
1500	75\\
1520	75\\
1540	75\\
1560	75\\
1580	75\\
1600	75\\
1620	75\\
1640	75\\
1660	75\\
1680	75\\
1700	75\\
1720	75\\
1740	75\\
1760	75\\
1780	75\\
1800	75\\
1820	75\\
1840	75\\
1860	75\\
1880	75\\
1900	75\\
1920	75\\
1940	75\\
1960	75\\
1980	75\\
2000	75\\
2020	75\\
2040	75\\
2060	75\\
2080	75\\
2100	75\\
2120	75\\
2140	75\\
2160	75\\
2180	75\\
2200	75\\
2220	75\\
2240	75\\
2260	75\\
2280	75\\
2300	75\\
2320	75\\
2340	75\\
2360	75\\
2380	75\\
2400	75\\
2420	75\\
2440	75\\
2460	75\\
2480	75\\
2500	75\\
2520	75\\
2540	75\\
2560	75\\
2580	75\\
2600	75\\
2620	75\\
2640	75\\
2660	75\\
2680	75\\
2700	75\\
2720	75\\
2740	75\\
2760	75\\
2780	75\\
2800	75\\
2820	75\\
2840	75\\
2860	75\\
2880	75\\
2900	75\\
2920	75\\
2940	75\\
2960	75\\
2980	75\\
3000	75\\
3020	75\\
3040	75\\
3060	75\\
3080	75\\
3100	75\\
3120	75\\
3140	75\\
3160	75\\
3180	75\\
3200	75\\
3220	75\\
3240	75\\
3260	75\\
3280	75\\
3300	75\\
3320	75\\
3340	75\\
3360	75\\
3380	75\\
3400	75\\
3420	75\\
3440	75\\
3460	75\\
3480	75\\
3500	75\\
3520	75\\
3540	75\\
3560	75\\
3580	75\\
3600	75\\
3620	75\\
3640	75\\
3660	75\\
3680	75\\
3700	75\\
3720	75\\
3740	75\\
3760	75\\
3780	75\\
3800	75\\
3820	75\\
3840	75\\
3860	75\\
3880	75\\
3900	75\\
3920	75\\
3940	75\\
3960	75\\
3980	75\\
4000	75\\
4020	75\\
4040	75\\
4060	75\\
4080	75\\
4100	75\\
4120	75\\
4140	75\\
4160	75\\
4180	75\\
4200	75\\
4220	75\\
4240	75\\
4260	75\\
4280	75\\
4300	75\\
4320	75\\
4340	75\\
4360	75\\
4380	75\\
4400	75\\
4420	75\\
4440	75\\
4460	75\\
4480	75\\
4500	75\\
4520	75\\
4540	75\\
4560	75\\
4580	75\\
4600	75\\
4620	75\\
4640	75\\
4660	75\\
4680	75\\
4700	75\\
4720	75\\
4740	75\\
4760	75\\
4780	75\\
4800	75\\
4820	75\\
4840	75\\
4860	75\\
4880	75\\
4900	75\\
4920	75\\
4940	75\\
4960	75\\
4980	75\\
5000	75\\
5020	75\\
5040	75\\
5060	75\\
5080	75\\
5100	75\\
5120	75\\
5140	75\\
5160	75\\
5180	75\\
5200	75\\
5220	75\\
5240	75\\
5260	75\\
5280	75\\
5300	75\\
5320	75\\
5340	75\\
5360	75\\
5380	75\\
5400	75\\
5420	75\\
5440	75\\
5460	75\\
5480	75\\
5500	75\\
5520	75\\
5540	75\\
5560	75\\
5580	75\\
5600	75\\
5620	75\\
5640	75\\
5660	75\\
5680	75\\
5700	75\\
5720	75\\
5740	75\\
5760	75\\
5780	75\\
5800	75\\
5820	75\\
5840	75\\
5860	75\\
5880	75\\
5900	75\\
5920	75\\
5940	75\\
5960	75\\
5980	75\\
6000	75\\
6020	75\\
6040	75\\
6060	75\\
6080	75\\
6100	75\\
6120	75\\
6140	75\\
6160	75\\
6180	75\\
6200	75\\
6220	75\\
6240	75\\
6260	75\\
6280	75\\
6300	75\\
6320	75\\
6340	75\\
6360	75\\
6380	75\\
6400	75\\
6420	75\\
6440	75\\
6460	75\\
6480	75\\
6500	75\\
6520	75\\
6540	75\\
6560	75\\
6580	75\\
6600	75\\
6620	75\\
6640	75\\
6660	75\\
6680	75\\
6700	75\\
6720	75\\
6740	75\\
6760	75\\
6780	75\\
6800	75\\
6820	75\\
6840	75\\
6860	75\\
6880	75\\
6900	75\\
6920	75\\
6940	75\\
6960	75\\
6980	75\\
7000	75\\
7020	75\\
7040	75\\
7060	75\\
7080	75\\
7100	75\\
7120	75\\
7140	75\\
7160	75\\
7180	75\\
7200	75\\
7220	75\\
7240	75\\
7260	75\\
7280	75\\
7300	75\\
7320	75\\
7340	75\\
7360	75\\
7380	75\\
7400	75\\
7420	75\\
7440	75\\
7460	75\\
7480	75\\
7500	75\\
7520	75\\
7540	75\\
7560	75\\
7580	75\\
7600	75\\
7620	75\\
7640	75\\
7660	75\\
7680	75\\
7700	75\\
7720	75\\
7740	75\\
7760	75\\
7780	75\\
7800	75\\
7820	75\\
7840	75\\
7860	75\\
7880	75\\
7900	75\\
7920	75\\
7940	75\\
7960	75\\
7980	75\\
8000	75\\
8020	75\\
8040	75\\
8060	75\\
8080	75\\
8100	75\\
8120	75\\
8140	75\\
8160	75\\
8180	75\\
8200	75\\
8220	75\\
8240	75\\
8260	75\\
8280	75\\
8300	75\\
8320	75\\
8340	75\\
8360	75\\
8380	75\\
8400	75\\
8420	75\\
8440	75\\
8460	75\\
8480	75\\
8500	75\\
8520	75\\
8540	75\\
8560	75\\
8580	75\\
8600	75\\
8620	75\\
8640	75\\
};
%\addlegendentry{}

\end{axis}

\begin{axis}[%
width=1.603in,
height=1.1in,
at={(3.436in,2.406in)},
scale only axis,
xmin=0,
xmax=8640,
ymin=149,
ymax=153,
ylabel style={font=\color{white!15!black}},
ylabel={Mass Flow @ Demand [kg/s]},
axis background/.style={fill=white},
legend style={legend cell align=left, align=left, draw=white!15!black}
]
\addplot [color=blue, line width=1.0pt]
  table[row sep=crcr]{%
0	150\\
20	150\\
40	150\\
60	150\\
80	150\\
100	150\\
120	150\\
140	150\\
160	150\\
180	150\\
200	150\\
220	150\\
240	150\\
260	150\\
280	150\\
300	150\\
320	150\\
340	150\\
360	150\\
380	150\\
400	150\\
420	150\\
440	150\\
460	150\\
480	150\\
500	150\\
520	150\\
540	150\\
560	150\\
580	150\\
600	150\\
620	150\\
640	150\\
660	150\\
680	150\\
700	150\\
720	150\\
740	150\\
760	150\\
780	150\\
800	150\\
820	150\\
840	150\\
860	150\\
880	150\\
900	150\\
920	150\\
940	150\\
960	151\\
980	151\\
1000	151\\
1020	151\\
1040	151\\
1060	151\\
1080	151\\
1100	151\\
1120	151\\
1140	151\\
1160	151\\
1180	151\\
1200	151\\
1220	151\\
1240	151\\
1260	151\\
1280	151\\
1300	151\\
1320	151\\
1340	151\\
1360	151\\
1380	151\\
1400	151\\
1420	151\\
1440	151\\
1460	151\\
1480	151\\
1500	151\\
1520	151\\
1540	151\\
1560	151\\
1580	151\\
1600	151\\
1620	151\\
1640	151\\
1660	151\\
1680	151\\
1700	151\\
1720	151\\
1740	151\\
1760	151\\
1780	151\\
1800	151\\
1820	151\\
1840	151\\
1860	151\\
1880	151\\
1900	151\\
1920	152\\
1940	152\\
1960	152\\
1980	152\\
2000	152\\
2020	152\\
2040	152\\
2060	152\\
2080	152\\
2100	152\\
2120	152\\
2140	152\\
2160	152\\
2180	152\\
2200	152\\
2220	152\\
2240	152\\
2260	152\\
2280	152\\
2300	152\\
2320	152\\
2340	152\\
2360	152\\
2380	152\\
2400	152\\
2420	152\\
2440	152\\
2460	152\\
2480	152\\
2500	152\\
2520	152\\
2540	152\\
2560	152\\
2580	152\\
2600	152\\
2620	152\\
2640	152\\
2660	152\\
2680	152\\
2700	152\\
2720	152\\
2740	152\\
2760	152\\
2780	152\\
2800	152\\
2820	152\\
2840	152\\
2860	152\\
2880	153\\
2900	153\\
2920	153\\
2940	153\\
2960	153\\
2980	153\\
3000	153\\
3020	153\\
3040	153\\
3060	153\\
3080	153\\
3100	153\\
3120	153\\
3140	153\\
3160	153\\
3180	153\\
3200	153\\
3220	153\\
3240	153\\
3260	153\\
3280	153\\
3300	153\\
3320	153\\
3340	153\\
3360	153\\
3380	153\\
3400	153\\
3420	153\\
3440	153\\
3460	153\\
3480	153\\
3500	153\\
3520	153\\
3540	153\\
3560	153\\
3580	153\\
3600	153\\
3620	153\\
3640	153\\
3660	153\\
3680	153\\
3700	153\\
3720	153\\
3740	153\\
3760	153\\
3780	153\\
3800	153\\
3820	153\\
3840	152\\
3860	152\\
3880	152\\
3900	152\\
3920	152\\
3940	152\\
3960	152\\
3980	152\\
4000	152\\
4020	152\\
4040	152\\
4060	152\\
4080	152\\
4100	152\\
4120	152\\
4140	152\\
4160	152\\
4180	152\\
4200	152\\
4220	152\\
4240	152\\
4260	152\\
4280	152\\
4300	152\\
4320	152\\
4340	152\\
4360	152\\
4380	152\\
4400	152\\
4420	152\\
4440	152\\
4460	152\\
4480	152\\
4500	152\\
4520	152\\
4540	152\\
4560	152\\
4580	152\\
4600	152\\
4620	152\\
4640	152\\
4660	152\\
4680	152\\
4700	152\\
4720	152\\
4740	152\\
4760	152\\
4780	152\\
4800	151\\
4820	151\\
4840	151\\
4860	151\\
4880	151\\
4900	151\\
4920	151\\
4940	151\\
4960	151\\
4980	151\\
5000	151\\
5020	151\\
5040	151\\
5060	151\\
5080	151\\
5100	151\\
5120	151\\
5140	151\\
5160	151\\
5180	151\\
5200	151\\
5220	151\\
5240	151\\
5260	151\\
5280	151\\
5300	151\\
5320	151\\
5340	151\\
5360	151\\
5380	151\\
5400	151\\
5420	151\\
5440	151\\
5460	151\\
5480	151\\
5500	151\\
5520	151\\
5540	151\\
5560	151\\
5580	151\\
5600	151\\
5620	151\\
5640	151\\
5660	151\\
5680	151\\
5700	151\\
5720	151\\
5740	151\\
5760	150\\
5780	150\\
5800	150\\
5820	150\\
5840	150\\
5860	150\\
5880	150\\
5900	150\\
5920	150\\
5940	150\\
5960	150\\
5980	150\\
6000	150\\
6020	150\\
6040	150\\
6060	150\\
6080	150\\
6100	150\\
6120	150\\
6140	150\\
6160	150\\
6180	150\\
6200	150\\
6220	150\\
6240	150\\
6260	150\\
6280	150\\
6300	150\\
6320	150\\
6340	150\\
6360	150\\
6380	150\\
6400	150\\
6420	150\\
6440	150\\
6460	150\\
6480	150\\
6500	150\\
6520	150\\
6540	150\\
6560	150\\
6580	150\\
6600	150\\
6620	150\\
6640	150\\
6660	150\\
6680	150\\
6700	150\\
6720	149\\
6740	149\\
6760	149\\
6780	149\\
6800	149\\
6820	149\\
6840	149\\
6860	149\\
6880	149\\
6900	149\\
6920	149\\
6940	149\\
6960	149\\
6980	149\\
7000	149\\
7020	149\\
7040	149\\
7060	149\\
7080	149\\
7100	149\\
7120	149\\
7140	149\\
7160	149\\
7180	149\\
7200	149\\
7220	149\\
7240	149\\
7260	149\\
7280	149\\
7300	149\\
7320	149\\
7340	149\\
7360	149\\
7380	149\\
7400	149\\
7420	149\\
7440	149\\
7460	149\\
7480	149\\
7500	149\\
7520	149\\
7540	149\\
7560	149\\
7580	149\\
7600	149\\
7620	149\\
7640	149\\
7660	149\\
7680	150\\
7700	150\\
7720	150\\
7740	150\\
7760	150\\
7780	150\\
7800	150\\
7820	150\\
7840	150\\
7860	150\\
7880	150\\
7900	150\\
7920	150\\
7940	150\\
7960	150\\
7980	150\\
8000	150\\
8020	150\\
8040	150\\
8060	150\\
8080	150\\
8100	150\\
8120	150\\
8140	150\\
8160	150\\
8180	150\\
8200	150\\
8220	150\\
8240	150\\
8260	150\\
8280	150\\
8300	150\\
8320	150\\
8340	150\\
8360	150\\
8380	150\\
8400	150\\
8420	150\\
8440	150\\
8460	150\\
8480	150\\
8500	150\\
8520	150\\
8540	150\\
8560	150\\
8580	150\\
8600	150\\
8620	150\\
8640	150\\
};
%\addlegendentry{}

\end{axis}

\begin{axis}[%
width=1.6in,
height=1.1in,
at={(1.011in,0.642in)},
scale only axis,
xmin=0,
xmax=8640,
ymin=148,
ymax=154,
ylabel style={font=\color{white!15!black}},
ylabel={Mass Flow @ Supply [kg/s]},
axis background/.style={fill=white},
legend style={legend cell align=left, align=left, draw=white!15!black}
]
\addplot [color=blue, line width=1.0pt]
  table[row sep=crcr]{%
0	150\\
20	150\\
40	150\\
60	150\\
80	150\\
100	150\\
120	150\\
140	150\\
160	150\\
180	150\\
200	150\\
220	150\\
240	150\\
260	150\\
280	150\\
300	150\\
320	150\\
340	150\\
360	150\\
380	150\\
400	150\\
420	150\\
440	150\\
460	150\\
480	150\\
500	150\\
520	150\\
540	150\\
560	150\\
580	150\\
600	150\\
620	150\\
640	150\\
660	150\\
680	150\\
700	150\\
720	150\\
740	150\\
760	150\\
780	150\\
800	150\\
820	150\\
840	150\\
860	150\\
880	150\\
900	150\\
920	150\\
940	150\\
960	150\\
980	151.010314595145\\
1000	150.999660250325\\
1020	150.999844796229\\
1040	150.999999991819\\
1060	150.999999999823\\
1080	150.999999999927\\
1100	150.999999999969\\
1120	151.00000000001\\
1140	151.000000000052\\
1160	151.000000000094\\
1180	151.000000000103\\
1200	151.00000000008\\
1220	151.000000000058\\
1240	151.000000000036\\
1260	151.000000000013\\
1280	150.999999999997\\
1300	150.999999999997\\
1320	150.999999999998\\
1340	150.999999999999\\
1360	150.999999999999\\
1380	151\\
1400	151\\
1420	151\\
1440	151\\
1460	151\\
1480	151\\
1500	151\\
1520	151\\
1540	151\\
1560	151\\
1580	151\\
1600	151\\
1620	151\\
1640	151\\
1660	151\\
1680	151\\
1700	151\\
1720	151\\
1740	151\\
1760	151\\
1780	151\\
1800	151\\
1820	151\\
1840	151\\
1860	151\\
1880	151\\
1900	151\\
1920	151\\
1940	152.010640363087\\
1960	151.999937355459\\
1980	151.999996467822\\
2000	152.000000042133\\
2020	152.00000007741\\
2040	152.000000048135\\
2060	152.000000039867\\
2080	152.0000000316\\
2100	152.000000023332\\
2120	152.000000015065\\
2140	152.000000006797\\
2160	151.99999999853\\
2180	151.999999991372\\
2200	151.99999999246\\
2220	151.999999993547\\
2240	151.999999994634\\
2260	151.999999995721\\
2280	151.999999996808\\
2300	151.999999997896\\
2320	151.999999998983\\
2340	151.999999999716\\
2360	151.999999999753\\
2380	151.99999999979\\
2400	151.999999999828\\
2420	151.999999999865\\
2440	151.999999999902\\
2460	151.999999999939\\
2480	151.999999999977\\
2500	151.999999999995\\
2520	151.999999999999\\
2540	152.000000000002\\
2560	152.000000000006\\
2580	152.000000000008\\
2600	152.000000000006\\
2620	152.000000000004\\
2640	152.000000000002\\
2660	152\\
2680	152\\
2700	152\\
2720	152\\
2740	152\\
2760	152\\
2780	152\\
2800	152\\
2820	152\\
2840	152\\
2860	152\\
2880	152\\
2900	153.01094167677\\
2920	152.99980343204\\
2940	152.99999760639\\
2960	152.999999563542\\
2980	152.999999936042\\
3000	153.000000176197\\
3020	153.000000099809\\
3040	153.00000002342\\
3060	153.000000005482\\
3080	153.000000002286\\
3100	153.000000000054\\
3120	153.000000000043\\
3140	153.000000000032\\
3160	153.000000000021\\
3180	153.00000000001\\
3200	152.999999999999\\
3220	152.999999999988\\
3240	152.999999999977\\
3260	152.999999999966\\
3280	152.999999999956\\
3300	152.999999999945\\
3320	152.999999999934\\
3340	152.999999999923\\
3360	152.999999999912\\
3380	152.999999999901\\
3400	152.99999999989\\
3420	152.999999999879\\
3440	152.999999999868\\
3460	152.999999999857\\
3480	152.999999999847\\
3500	152.999999999836\\
3520	152.999999999825\\
3540	152.999999999814\\
3560	152.999999999803\\
3580	152.999999999792\\
3600	152.999999999812\\
3620	152.999999999893\\
3640	152.999999999974\\
3660	152.999999999992\\
3680	152.999999999998\\
3700	153.000000000002\\
3720	153.000000000001\\
3740	153\\
3760	153\\
3780	153\\
3800	153\\
3820	153\\
3840	153\\
3860	151.989139582213\\
3880	152.000252623162\\
3900	151.999999750886\\
3920	151.999999901515\\
3940	151.999999984245\\
3960	151.999999993372\\
3980	152.0000000025\\
4000	152.000000011628\\
4020	152.000000008626\\
4040	152.000000005577\\
4060	152.000000002529\\
4080	151.999999999876\\
4100	151.999999999901\\
4120	151.999999999927\\
4140	151.999999999952\\
4160	151.999999999973\\
4180	151.999999999979\\
4200	151.999999999985\\
4220	151.999999999991\\
4240	151.999999999996\\
4260	151.999999999997\\
4280	151.999999999998\\
4300	152\\
4320	152.000000000001\\
4340	152\\
4360	152\\
4380	152\\
4400	152\\
4420	152\\
4440	152\\
4460	152\\
4480	152\\
4500	152\\
4520	152\\
4540	152\\
4560	152\\
4580	152\\
4600	152\\
4620	152\\
4640	152\\
4660	152\\
4680	152\\
4700	152\\
4720	152\\
4740	152\\
4760	152\\
4780	152\\
4800	152\\
4820	150.989472982865\\
4840	151.000385568715\\
4860	151.000376530526\\
4880	151.000000012996\\
4900	150.999999996487\\
4920	150.999999999799\\
4940	150.99999999986\\
4960	150.999999999921\\
4980	150.999999999981\\
5000	151.000000000042\\
5020	151.00000000005\\
5040	151.000000000039\\
5060	151.000000000028\\
5080	151.000000000016\\
5100	151.000000000005\\
5120	151.000000000002\\
5140	151.000000000002\\
5160	151.000000000001\\
5180	151.000000000001\\
5200	151\\
5220	151\\
5240	151\\
5260	151\\
5280	151\\
5300	151\\
5320	151\\
5340	151\\
5360	151\\
5380	151\\
5400	151\\
5420	151\\
5440	151\\
5460	151\\
5480	151\\
5500	151\\
5520	151\\
5540	151\\
5560	151\\
5580	151\\
5600	151\\
5620	151\\
5640	151\\
5660	151\\
5680	151\\
5700	151\\
5720	151\\
5740	151\\
5760	151\\
5780	149.989784707729\\
5800	150.000216845846\\
5820	150.000237825012\\
5840	149.999999919024\\
5860	150.000000009327\\
5880	150.000000003157\\
5900	149.999999998475\\
5920	149.999999999119\\
5940	149.999999999763\\
5960	149.999999999933\\
5980	149.999999999982\\
6000	150.000000000012\\
6020	150.000000000012\\
6040	150.000000000011\\
6060	150.000000000011\\
6080	150.00000000001\\
6100	150.00000000001\\
6120	150.000000000009\\
6140	150.000000000009\\
6160	150.000000000008\\
6180	150.000000000008\\
6200	150.000000000007\\
6220	150.000000000007\\
6240	150.000000000006\\
6260	150.000000000006\\
6280	150.000000000005\\
6300	150.000000000005\\
6320	150.000000000004\\
6340	150.000000000003\\
6360	150.000000000003\\
6380	150.000000000003\\
6400	150.000000000002\\
6420	150.000000000001\\
6440	150.000000000001\\
6460	150\\
6480	150\\
6500	150\\
6520	150\\
6540	150\\
6560	150\\
6580	150\\
6600	150\\
6620	150\\
6640	150\\
6660	150\\
6680	150\\
6700	150\\
6720	150\\
6740	148.99014439314\\
6760	149.000133252551\\
6780	149.000001609047\\
6800	148.999999738765\\
6820	149.000000023868\\
6840	149.000000018453\\
6860	149.000000006865\\
6880	149.000000000596\\
6900	149.00000000053\\
6920	149.000000000465\\
6940	149.000000000399\\
6960	149.000000000334\\
6980	149.000000000269\\
7000	149.000000000203\\
7020	149.000000000138\\
7040	149.000000000072\\
7060	149.000000000007\\
7080	148.999999999941\\
7100	148.999999999876\\
7120	148.99999999981\\
7140	148.999999999745\\
7160	148.999999999724\\
7180	148.999999999745\\
7200	148.999999999765\\
7220	148.999999999786\\
7240	148.999999999807\\
7260	148.999999999827\\
7280	148.999999999848\\
7300	148.999999999868\\
7320	148.999999999889\\
7340	148.99999999991\\
7360	148.99999999993\\
7380	148.999999999951\\
7400	148.999999999972\\
7420	148.999999999992\\
7440	149.000000000002\\
7460	149.000000000002\\
7480	149\\
7500	149\\
7520	149\\
7540	149\\
7560	149\\
7580	149\\
7600	149\\
7620	149\\
7640	149\\
7660	149\\
7680	149\\
7700	150.009984643847\\
7720	149.99978599258\\
7740	149.99978382006\\
7760	150.000000401546\\
7780	150.000000087374\\
7800	149.999999995648\\
7820	150.000000000126\\
7840	150.00000000011\\
7860	150.000000000094\\
7880	150.000000000078\\
7900	150.000000000062\\
7920	150.000000000047\\
7940	150.000000000031\\
7960	150.000000000025\\
7980	150.00000000002\\
8000	150.000000000016\\
8020	150.000000000012\\
8040	150.000000000008\\
8060	150.000000000003\\
8080	150\\
8100	150\\
8120	150\\
8140	150\\
8160	150\\
8180	150\\
8200	150\\
8220	150\\
8240	150\\
8260	150\\
8280	150\\
8300	150\\
8320	150\\
8340	150\\
8360	150\\
8380	150\\
8400	150\\
8420	150\\
8440	150\\
8460	150\\
8480	150\\
8500	150\\
8520	150\\
8540	150\\
8560	150\\
8580	150\\
8600	150\\
8620	150\\
8640	150\\
};
%\addlegendentry{}

\end{axis}

\begin{axis}[%
width=1.603in,
height=1.1in,
at={(3.436in,0.642in)},
scale only axis,
xmin=0,
xmax=8640,
ymin=72.55,
ymax=72.7069954426378,
ylabel style={font=\color{white!15!black}},
ylabel={Pressure @ Demand [bar]},
axis background/.style={fill=white},
legend style={legend cell align=left, align=left, draw=white!15!black}
]
\addplot [color=red, line width=1.0pt]
  table[row sep=crcr]{%
0	72.6756179326083\\
20	72.6756179326083\\
40	72.6756179326083\\
60	72.6756179326083\\
80	72.6756179326083\\
100	72.6756179326083\\
120	72.6756179326083\\
140	72.6756179326083\\
160	72.6756179326083\\
180	72.6756179326083\\
200	72.6756179326083\\
220	72.6756179326083\\
240	72.6756179326083\\
260	72.6756179326083\\
280	72.6756179326083\\
300	72.6756179326083\\
320	72.6756179326083\\
340	72.6756179326083\\
360	72.6756179326083\\
380	72.6756179326083\\
400	72.6756179326083\\
420	72.6756179326083\\
440	72.6756179326083\\
460	72.6756179326083\\
480	72.6756179326083\\
500	72.6756179326083\\
520	72.6756179326083\\
540	72.6756179326083\\
560	72.6756179326083\\
580	72.6756179326083\\
600	72.6756179326083\\
620	72.6756179326083\\
640	72.6756179326083\\
660	72.6756179326083\\
680	72.6756179326083\\
700	72.6756179326083\\
720	72.6756179326083\\
740	72.6756179326083\\
760	72.6756179326083\\
780	72.6756179326083\\
800	72.6756179326083\\
820	72.6756179326083\\
840	72.6756179326083\\
860	72.6756179326083\\
880	72.6756179326083\\
900	72.6756179326083\\
920	72.6756179326083\\
940	72.6756179326083\\
960	72.6755469782406\\
980	72.6441353997748\\
1000	72.6440152633452\\
1020	72.6440167744887\\
1040	72.6440186410467\\
1060	72.6440186404945\\
1080	72.6440186404665\\
1100	72.6440186404676\\
1120	72.6440186404687\\
1140	72.6440186404698\\
1160	72.6440186404709\\
1180	72.6440186404715\\
1200	72.6440186404717\\
1220	72.6440186404719\\
1240	72.6440186404721\\
1260	72.6440186404723\\
1280	72.6440186404724\\
1300	72.6440186404724\\
1320	72.6440186404724\\
1340	72.6440186404723\\
1360	72.6440186404723\\
1380	72.6440186404723\\
1400	72.6440186404723\\
1420	72.6440186404723\\
1440	72.6440186404723\\
1460	72.6440186404723\\
1480	72.6440186404723\\
1500	72.6440186404723\\
1520	72.6440186404723\\
1540	72.6440186404723\\
1560	72.6440186404723\\
1580	72.6440186404723\\
1600	72.6440186404723\\
1620	72.6440186404723\\
1640	72.6440186404723\\
1660	72.6440186404723\\
1680	72.6440186404723\\
1700	72.6440186404723\\
1720	72.6440186404723\\
1740	72.6440186404723\\
1760	72.6440186404723\\
1780	72.6440186404723\\
1800	72.6440186404723\\
1820	72.6440186404723\\
1840	72.6440186404723\\
1860	72.6440186404723\\
1880	72.6440186404723\\
1900	72.6440186404723\\
1920	72.6439912802173\\
1940	72.6123009229617\\
1960	72.6121936948696\\
1980	72.612195529005\\
2000	72.6121954855004\\
2020	72.612195490208\\
2040	72.6121954907433\\
2060	72.612195490808\\
2080	72.6121954908728\\
2100	72.6121954909375\\
2120	72.6121954910022\\
2140	72.612195491067\\
2160	72.6121954911317\\
2180	72.6121954911875\\
2200	72.6121954911776\\
2220	72.6121954911677\\
2240	72.6121954911577\\
2260	72.6121954911477\\
2280	72.6121954911378\\
2300	72.6121954911278\\
2320	72.6121954911179\\
2340	72.6121954911104\\
2360	72.6121954911079\\
2380	72.6121954911053\\
2400	72.6121954911028\\
2420	72.6121954911002\\
2440	72.6121954910977\\
2460	72.6121954910951\\
2480	72.6121954910925\\
2500	72.6121954910916\\
2520	72.6121954910919\\
2540	72.6121954910923\\
2560	72.6121954910926\\
2580	72.6121954910929\\
2600	72.6121954910929\\
2620	72.612195491093\\
2640	72.612195491093\\
2660	72.6121954910931\\
2680	72.6121954910931\\
2700	72.6121954910931\\
2720	72.6121954910931\\
2740	72.6121954910931\\
2760	72.6121954910931\\
2780	72.6121954910931\\
2800	72.6121954910931\\
2820	72.6121954910931\\
2840	72.6121954910931\\
2860	72.6121954910931\\
2880	72.6120976909684\\
2900	72.580242460277\\
2920	72.58014569711\\
2940	72.5801482053118\\
2960	72.5801481873805\\
2980	72.5801481879015\\
3000	72.5801481884115\\
3020	72.5801481888953\\
3040	72.580148189379\\
3060	72.580148189561\\
3080	72.5801481896668\\
3100	72.5801481897401\\
3120	72.5801481897384\\
3140	72.5801481897367\\
3160	72.5801481897351\\
3180	72.5801481897334\\
3200	72.5801481897317\\
3220	72.58014818973\\
3240	72.5801481897283\\
3260	72.5801481897266\\
3280	72.5801481897249\\
3300	72.5801481897232\\
3320	72.5801481897215\\
3340	72.5801481897199\\
3360	72.5801481897182\\
3380	72.5801481897165\\
3400	72.5801481897148\\
3420	72.5801481897131\\
3440	72.5801481897114\\
3460	72.5801481897097\\
3480	72.580148189708\\
3500	72.5801481897063\\
3520	72.5801481897047\\
3540	72.580148189703\\
3560	72.5801481897013\\
3580	72.5801481896996\\
3600	72.5801481896981\\
3620	72.5801481896968\\
3640	72.5801481896956\\
3660	72.5801481896954\\
3680	72.5801481896954\\
3700	72.5801481896954\\
3720	72.5801481896954\\
3740	72.5801481896954\\
3760	72.5801481896954\\
3780	72.5801481896954\\
3800	72.5801481896954\\
3820	72.5801481896954\\
3840	72.5803065428281\\
3860	72.6120936343938\\
3880	72.6121976892357\\
3900	72.6121955321089\\
3920	72.6121954908908\\
3940	72.6121954906516\\
3960	72.6121954907914\\
3980	72.6121954909313\\
4000	72.6121954910712\\
4020	72.6121954910783\\
4040	72.6121954910848\\
4060	72.6121954910913\\
4080	72.6121954910969\\
4100	72.6121954910964\\
4120	72.612195491096\\
4140	72.6121954910955\\
4160	72.612195491095\\
4180	72.6121954910945\\
4200	72.612195491094\\
4220	72.6121954910935\\
4240	72.6121954910931\\
4260	72.6121954910931\\
4280	72.6121954910931\\
4300	72.6121954910931\\
4320	72.612195491093\\
4340	72.6121954910931\\
4360	72.6121954910931\\
4380	72.6121954910931\\
4400	72.6121954910931\\
4420	72.6121954910931\\
4440	72.6121954910931\\
4460	72.6121954910931\\
4480	72.6121954910931\\
4500	72.6121954910931\\
4520	72.6121954910931\\
4540	72.6121954910931\\
4560	72.6121954910931\\
4580	72.6121954910931\\
4600	72.6121954910931\\
4620	72.6121954910931\\
4640	72.6121954910931\\
4660	72.6121954910931\\
4680	72.6121954910931\\
4700	72.6121954910931\\
4720	72.6121954910931\\
4740	72.6121954910931\\
4760	72.6121954910931\\
4780	72.6121954910931\\
4800	72.6122260030897\\
4820	72.6439055006662\\
4840	72.6440223148718\\
4860	72.6440229743784\\
4880	72.6440186387449\\
4900	72.6440186404516\\
4920	72.6440186404702\\
4940	72.6440186404708\\
4960	72.6440186404714\\
4980	72.6440186404721\\
5000	72.6440186404727\\
5020	72.6440186404727\\
5040	72.6440186404726\\
5060	72.6440186404725\\
5080	72.6440186404724\\
5100	72.6440186404722\\
5120	72.6440186404722\\
5140	72.6440186404722\\
5160	72.6440186404723\\
5180	72.6440186404723\\
5200	72.6440186404723\\
5220	72.6440186404723\\
5240	72.6440186404723\\
5260	72.6440186404723\\
5280	72.6440186404723\\
5300	72.6440186404723\\
5320	72.6440186404723\\
5340	72.6440186404723\\
5360	72.6440186404723\\
5380	72.6440186404723\\
5400	72.6440186404723\\
5420	72.6440186404723\\
5440	72.6440186404723\\
5460	72.6440186404723\\
5480	72.6440186404723\\
5500	72.6440186404723\\
5520	72.6440186404723\\
5540	72.6440186404723\\
5560	72.6440186404723\\
5580	72.6440186404723\\
5600	72.6440186404723\\
5620	72.6440186404723\\
5640	72.6440186404723\\
5660	72.6440186404723\\
5680	72.6440186404723\\
5700	72.6440186404723\\
5720	72.6440186404723\\
5740	72.6440186404723\\
5760	72.6440995286175\\
5780	72.6754930565247\\
5800	72.6756202995842\\
5820	72.6756200008125\\
5840	72.6756179311862\\
5860	72.6756179326834\\
5880	72.6756179326463\\
5900	72.6756179326164\\
5920	72.6756179326124\\
5940	72.6756179326083\\
5960	72.6756179326077\\
5980	72.675617932608\\
6000	72.6756179326082\\
6020	72.6756179326082\\
6040	72.6756179326082\\
6060	72.6756179326082\\
6080	72.6756179326082\\
6100	72.6756179326082\\
6120	72.6756179326082\\
6140	72.6756179326083\\
6160	72.6756179326083\\
6180	72.6756179326083\\
6200	72.6756179326083\\
6220	72.6756179326083\\
6240	72.6756179326083\\
6260	72.6756179326083\\
6280	72.6756179326083\\
6300	72.6756179326083\\
6320	72.6756179326083\\
6340	72.6756179326083\\
6360	72.6756179326083\\
6380	72.6756179326083\\
6400	72.6756179326083\\
6420	72.6756179326083\\
6440	72.6756179326083\\
6460	72.6756179326083\\
6480	72.6756179326083\\
6500	72.6756179326083\\
6520	72.6756179326083\\
6540	72.6756179326083\\
6560	72.6756179326083\\
6580	72.6756179326083\\
6600	72.6756179326083\\
6620	72.6756179326083\\
6640	72.6756179326083\\
6660	72.6756179326083\\
6680	72.6756179326083\\
6700	72.6756179326083\\
6720	72.6756225089805\\
6740	72.7068570978263\\
6760	72.7069954426378\\
6780	72.7069936959585\\
6800	72.7069936625635\\
6820	72.7069936589782\\
6840	72.7069936594826\\
6860	72.7069936596472\\
6880	72.7069936597088\\
6900	72.706993659707\\
6920	72.7069936597052\\
6940	72.7069936597034\\
6960	72.7069936597016\\
6980	72.7069936596999\\
7000	72.7069936596981\\
7020	72.7069936596963\\
7040	72.7069936596945\\
7060	72.7069936596927\\
7080	72.7069936596909\\
7100	72.7069936596891\\
7120	72.7069936596874\\
7140	72.7069936596856\\
7160	72.7069936596847\\
7180	72.7069936596846\\
7200	72.7069936596846\\
7220	72.7069936596846\\
7240	72.7069936596845\\
7260	72.7069936596845\\
7280	72.7069936596844\\
7300	72.7069936596844\\
7320	72.7069936596843\\
7340	72.7069936596843\\
7360	72.7069936596843\\
7380	72.7069936596842\\
7400	72.7069936596842\\
7420	72.7069936596841\\
7440	72.7069936596842\\
7460	72.7069936596843\\
7480	72.7069936596843\\
7500	72.7069936596844\\
7520	72.7069936596844\\
7540	72.7069936596844\\
7560	72.7069936596844\\
7580	72.7069936596844\\
7600	72.7069936596844\\
7620	72.7069936596844\\
7640	72.7069936596844\\
7660	72.7069936596844\\
7680	72.706857830168\\
7700	72.6757462947022\\
7720	72.6756157273676\\
7740	72.6756163181801\\
7760	72.6756179302789\\
7780	72.6756179333653\\
7800	72.6756179325681\\
7820	72.6756179326085\\
7840	72.6756179326085\\
7860	72.6756179326085\\
7880	72.6756179326085\\
7900	72.6756179326085\\
7920	72.6756179326086\\
7940	72.6756179326086\\
7960	72.6756179326086\\
7980	72.6756179326085\\
8000	72.6756179326085\\
8020	72.6756179326084\\
8040	72.6756179326084\\
8060	72.6756179326083\\
8080	72.6756179326083\\
8100	72.6756179326083\\
8120	72.6756179326083\\
8140	72.6756179326083\\
8160	72.6756179326083\\
8180	72.6756179326083\\
8200	72.6756179326083\\
8220	72.6756179326083\\
8240	72.6756179326083\\
8260	72.6756179326083\\
8280	72.6756179326083\\
8300	72.6756179326083\\
8320	72.6756179326083\\
8340	72.6756179326083\\
8360	72.6756179326083\\
8380	72.6756179326083\\
8400	72.6756179326083\\
8420	72.6756179326083\\
8440	72.6756179326083\\
8460	72.6756179326083\\
8480	72.6756179326083\\
8500	72.6756179326083\\
8520	72.6756179326083\\
8540	72.6756179326083\\
8560	72.6756179326083\\
8580	72.6756179326083\\
8600	72.6756179326083\\
8620	72.6756179326083\\
8640	72.6756179326083\\
};
%\addlegendentry{flow}

\end{axis}
\end{tikzpicture}%

%% file: hyp_diamond_step.tex
\definecolor{rasp}{rgb}{.89,.04,.36}      % rot fuer mid
\begin{tikzpicture}

\begin{axis}[%
width=1.603in,
height=1.1in,
at={(1.011in,2.406in)},
scale only axis,
xmin=0,
xmax=100,
ymin=69,
ymax=71,
ylabel style={font=\color{white!15!black}},
ylabel={Pressure @ Supply [bar]},
axis background/.style={fill=white},
legend style={legend cell align=left, align=left, draw=white!15!black}
]
\addplot [color=red, line width=1.0pt]
  table[row sep=crcr]{%
0	70\\
1	70\\
2	70\\
3	70\\
4	70\\
5	70\\
6	70\\
7	70\\
8	70\\
9	70\\
10	70\\
11	70\\
12	70\\
13	70\\
14	70\\
15	70\\
16	70\\
17	70\\
18	70\\
19	70\\
20	70\\
21	70\\
22	70\\
23	70\\
24	70\\
25	70\\
26	70\\
27	70\\
28	70\\
29	70\\
30	70\\
31	70\\
32	70\\
33	70\\
34	70\\
35	70\\
36	70\\
37	70\\
38	70\\
39	70\\
40	70\\
41	70\\
42	70\\
43	70\\
44	70\\
45	70\\
46	70\\
47	70\\
48	70\\
49	70\\
50	70\\
51	70\\
52	70\\
53	70\\
54	70\\
55	70\\
56	70\\
57	70\\
58	70\\
59	70\\
60	70\\
61	70\\
62	70\\
63	70\\
64	70\\
65	70\\
66	70\\
67	70\\
68	70\\
69	70\\
70	70\\
71	70\\
72	70\\
73	70\\
74	70\\
75	70\\
76	70\\
77	70\\
78	70\\
79	70\\
80	70\\
81	70\\
82	70\\
83	70\\
84	70\\
85	70\\
86	70\\
87	70\\
88	70\\
89	70\\
90	70\\
91	70\\
92	70\\
93	70\\
94	70\\
95	70\\
96	70\\
97	70\\
98	70\\
99	70\\
100	70\\
};
%\addlegendentry{end}

\end{axis}

\begin{axis}[%
width=1.603in,
height=1.1in,
at={(3.236in,2.406in)},
scale only axis,
xmin=0,
xmax=100,
ymin=25,
ymax=50,
ylabel style={font=\color{white!15!black}},
ylabel={Mass Flow @ Demand [kg/s]},
axis background/.style={fill=white},
legend style={legend cell align=left, align=left, draw=white!15!black}
]
\addplot [color=blue, line width=1.0pt]
  table[row sep=crcr]{%
0	30\\
1	30\\
2	30\\
3	30\\
4	30\\
5	30\\
6	30\\
7	30\\
8	30\\
9	30\\
10	30\\
11	30\\
12	30\\
13	30\\
14	30\\
15	30\\
16	30\\
17	30\\
18	30\\
19	30\\
20	40\\
21	40\\
22	40\\
23	40\\
24	40\\
25	40\\
26	40\\
27	40\\
28	40\\
29	40\\
30	40\\
31	40\\
32	40\\
33	40\\
34	40\\
35	40\\
36	40\\
37	40\\
38	40\\
39	40\\
40	40\\
41	40\\
42	40\\
43	40\\
44	40\\
45	40\\
46	40\\
47	40\\
48	40\\
49	40\\
50	40\\
51	40\\
52	40\\
53	40\\
54	40\\
55	40\\
56	40\\
57	40\\
58	40\\
59	40\\
60	40\\
61	40\\
62	40\\
63	40\\
64	40\\
65	40\\
66	40\\
67	40\\
68	40\\
69	40\\
70	40\\
71	40\\
72	40\\
73	40\\
74	40\\
75	40\\
76	40\\
77	40\\
78	40\\
79	40\\
80	40\\
81	40\\
82	40\\
83	40\\
84	40\\
85	40\\
86	40\\
87	40\\
88	40\\
89	40\\
90	40\\
91	40\\
92	40\\
93	40\\
94	40\\
95	40\\
96	40\\
97	40\\
98	40\\
99	40\\
100	40\\
};
%\addlegendentry{end}

\end{axis}

\begin{axis}[%
width=1.603in,
height=1.1in,
at={(1.011in,0.642in)},
scale only axis,
xmin=0,
xmax=100,
ymin=25,
ymax=50,
ylabel style={font=\color{white!15!black}},
ylabel={Mass Flow @ Supply [kg/s]},
axis background/.style={fill=white},
legend style={legend cell align=left, align=left, draw=white!15!black}
]

\addplot [color=rasp, line width=1.0pt]
  table[row sep=crcr]{%
0	29.9932120239495\\
1	29.993073792163\\
2	29.9939634558671\\
3	29.9956606807996\\
4	29.9972092259626\\
5	29.9984647568773\\
6	29.9989082899778\\
7	29.9997069330672\\
8	30.0007421532913\\
9	30.001586881257\\
10	30.002418311771\\
11	30.0031419255253\\
12	30.0044233593415\\
13	30.0048279468753\\
14	30.0042019958791\\
15	30.00311800682\\
16	30.0026889484648\\
17	30.0036775000516\\
18	30.0050723772706\\
19	30.0055085149278\\
20	30.0048913675557\\
21	30.1057966168413\\
22	29.9222417035133\\
23	30.630452806288\\
24	29.7363581728628\\
25	30.7555480710162\\
26	31.0437257174296\\
27	28.0354519278352\\
28	31.2337002640659\\
29	34.2430001892333\\
30	42.5090544314279\\
31	43.80198937205\\
32	42.5300100840341\\
33	41.3534761310263\\
34	43.4445081615521\\
35	44.8929492542119\\
36	45.1995969105088\\
37	44.56978554195\\
38	45.2189513251024\\
39	45.4228393773321\\
40	44.9128929367663\\
41	46.1080368961856\\
42	45.6743007144826\\
43	47.7465210784351\\
44	46.5731213643182\\
45	48.47552781216\\
46	47.3330310945388\\
47	45.8761916169268\\
48	41.1691492080749\\
49	39.8740017622139\\
50	39.3324267086339\\
51	41.4839300617027\\
52	40.406132907841\\
53	38.0675586261677\\
54	38.1826138574287\\
55	37.9182330400682\\
56	38.7065888591556\\
57	38.2515766871132\\
58	37.6821648450252\\
59	37.8722697134201\\
60	36.6978110706103\\
61	35.7727656096079\\
62	35.3734813778757\\
63	34.9046428141573\\
64	35.1493459252324\\
65	36.0414591046994\\
66	37.798806170049\\
67	39.2815617350704\\
68	39.626455259878\\
69	38.4463504735635\\
70	38.5997891859668\\
71	39.448376495822\\
72	41.0186886395273\\
73	40.2297045655888\\
74	39.9096972883478\\
75	39.7284284385969\\
76	40.7016516655908\\
77	40.5591805027014\\
78	41.4078651605248\\
79	41.7149065703731\\
80	43.1323603776246\\
81	42.5787804490832\\
82	42.8890597698242\\
83	42.4908791898219\\
84	41.4757026021274\\
85	41.1576814551837\\
86	40.5127453406579\\
87	41.4136338685345\\
88	41.4047815239476\\
89	41.1190634454667\\
90	39.9750117816991\\
91	40.2174390941634\\
92	40.4153151652202\\
93	41.1059927749889\\
94	40.3764824969941\\
95	39.9192306568402\\
96	39.7317048674447\\
97	39.3379954809661\\
98	38.7327541251992\\
99	38.5586297933143\\
100	38.5369988607852\\
};
\addlegendentry{mid}

\addplot [color=green, line width=1.0pt]
  table[row sep=crcr]{%
0	29.9923032307969\\
1	29.9919168792977\\
2	29.9914931670452\\
3	29.9928884999426\\
4	29.9943279197206\\
5	29.9962908639898\\
6	29.9961085783706\\
7	29.9954801685144\\
8	29.995347119443\\
9	29.9967710322603\\
10	29.9992877990816\\
11	30.0008828303685\\
12	30.0004331746184\\
13	29.9996033888453\\
14	29.9999538506655\\
15	30.0013869442782\\
16	30.0026350543671\\
17	30.0029918769005\\
18	30.0037992885234\\
19	30.0044845260354\\
20	30.005363544854\\
21	30.0051782961923\\
22	30.0042886153867\\
23	30.0035888897228\\
24	30.0030776824369\\
25	30.0039473761779\\
26	30.0039558046689\\
27	30.0088871067468\\
28	30.2411299715689\\
29	33.0448679856846\\
30	40.898932148083\\
31	43.6420404223224\\
32	45.6202438257567\\
33	40.6807868135703\\
34	43.520589684075\\
35	44.3544106599547\\
36	45.0702573650387\\
37	45.4781172843123\\
38	45.2765879714849\\
39	45.8602006071414\\
40	44.7767555302652\\
41	43.9151333397292\\
42	46.1472813413527\\
43	47.0948617265728\\
44	47.8209151764396\\
45	46.1500829581099\\
46	47.1692346311827\\
47	46.498009909956\\
48	44.513543294601\\
49	41.1320956417467\\
50	39.3678319852548\\
51	39.687059993099\\
52	40.1091105740466\\
53	39.9545939474744\\
54	39.5095666076242\\
55	37.3682040731341\\
56	37.3620144941663\\
57	37.7011638064017\\
58	38.3735452613579\\
59	38.3760890599663\\
60	38.0549356599968\\
61	37.1912113578585\\
62	35.4343151899348\\
63	35.0070065718342\\
64	35.505228680491\\
65	36.220165003006\\
66	36.7753782502951\\
67	36.0944272862818\\
68	38.7774244103171\\
69	38.9660031560232\\
70	39.9112640064206\\
71	38.5968056807239\\
72	38.9651430714438\\
73	40.4998091987351\\
74	40.8667246280106\\
75	40.6511816236993\\
76	40.672418860704\\
77	40.3694504281922\\
78	40.2610319255852\\
79	40.3599774126882\\
80	41.7002746093238\\
81	42.8023461899325\\
82	43.0017021972648\\
83	42.3400443854675\\
84	42.0074007563726\\
85	42.0657886199714\\
86	41.8411685204138\\
87	41.6317952783234\\
88	41.1237998988601\\
89	40.6900820018766\\
90	40.9724339957765\\
91	40.3989733216521\\
92	40.5904587791221\\
93	39.7164629488767\\
94	39.9951674953677\\
95	40.232480980375\\
96	40.3647421300043\\
97	40.2649837403095\\
98	39.7412860455625\\
99	39.4151868797641\\
100	38.6244346088101\\
};
\addlegendentry{end}
\addplot [color=black, line width=1.0pt]
  table[row sep=crcr]{%
0	30.0024102902499\\
1	29.9992238728098\\
2	29.9977386778601\\
3	30.0013616448288\\
4	30.0021852177677\\
5	29.9983177049597\\
6	29.9982311092462\\
7	30.0021726886998\\
8	30.001512831128\\
9	29.9976786084444\\
10	29.9990464785378\\
11	30.0025860814539\\
12	30.0005637159573\\
13	29.9974506626382\\
14	29.9999934969782\\
15	30.0025233326605\\
16	29.9995682655223\\
17	29.9976829547639\\
18	30.0008594143725\\
19	30.0020153569027\\
20	29.9987578965769\\
21	29.9983158649204\\
22	30.0014561875845\\
23	30.0011912525279\\
24	29.9983166760598\\
25	29.9995357335394\\
26	30.0109434637825\\
27	30.1160704578212\\
28	30.7673718924591\\
29	33.0119066353612\\
30	37.4409450639306\\
31	42.2714865840629\\
32	44.7202555194156\\
33	45.0500999477342\\
34	45.5052986347795\\
35	46.1036219777098\\
36	45.9333136611314\\
37	45.5833592862351\\
38	45.6024630273885\\
39	45.855927984514\\
40	45.9489652003429\\
41	45.6660964458512\\
42	45.8994264794145\\
43	46.622923913179\\
44	46.3647108495001\\
45	45.7100039520438\\
46	45.6902926877438\\
47	44.97486560166\\
48	43.2963987516569\\
49	41.9984786188136\\
50	41.0367692104815\\
51	40.0000122235945\\
52	39.035616626351\\
53	38.1906247022683\\
54	37.8604458380006\\
55	37.8492751376193\\
56	37.3622141080235\\
57	36.9439206523174\\
58	37.1331355363672\\
59	37.0672819634685\\
60	36.6111126705061\\
61	36.472396882563\\
62	36.5261266191101\\
63	36.5107347971601\\
64	36.6033812155272\\
65	36.8038790119277\\
66	37.162287986851\\
67	37.7677348161646\\
68	38.3236389881504\\
69	38.7132308850242\\
70	39.2405407733641\\
71	39.8785591144271\\
72	40.2834102786922\\
73	40.483097931485\\
74	40.7656394412435\\
75	41.1420251065737\\
76	41.3718142624411\\
77	41.3968908370004\\
78	41.5039614598667\\
79	41.8243893188416\\
80	41.9797266010998\\
81	41.8174094640267\\
82	41.8085836934498\\
83	42.0230572696088\\
84	41.9068897918545\\
85	41.4956915868697\\
86	41.3494571646085\\
87	41.3604514397548\\
88	41.0435306683472\\
89	40.5664822554906\\
90	40.3940544611679\\
91	40.3868111805439\\
92	40.1001541526203\\
93	39.6855441842482\\
94	39.6114406376279\\
95	39.6807998955039\\
96	39.4209351074477\\
97	39.1013650544128\\
98	39.1445305234906\\
99	39.2226594724639\\
100	39.0121501385661\\
};
\addlegendentry{new}
\end{axis}

\begin{axis}[%
width=1.603in,
height=1.1in,
at={(3.236in,0.642in)},
scale only axis,
xmin=0,
xmax=100,
ymin=69,
ymax=71,
ylabel style={font=\color{white!15!black}},
ylabel={Pressure @ Demand [bar]},
axis background/.style={fill=white},
legend style={legend cell align=left, align=left, draw=white!15!black}
]

 \addplot [color=rasp, line width=1.0pt]
  table[row sep=crcr]{%
0	69.9583683244541\\
1	69.9580098668402\\
2	69.9577818400878\\
3	69.957825072735\\
4	69.9581040504067\\
5	69.9584403626709\\
6	69.9586317937788\\
7	69.9585353533498\\
8	69.9582038695583\\
9	69.9578493560703\\
10	69.9577333749957\\
11	69.9579311706151\\
12	69.9582870334525\\
13	69.9585892254091\\
14	69.958613676222\\
15	69.9583630763739\\
16	69.9580073393235\\
17	69.9577747920053\\
18	69.9578488036611\\
19	69.9581653746835\\
20	69.958540599891\\
21	69.9423521581397\\
22	69.9087845122084\\
23	69.8406821696669\\
24	69.7577298828143\\
25	69.7918106836096\\
26	69.8208242385932\\
27	69.947350024261\\
28	69.9143044003184\\
29	69.8985460489424\\
30	69.7902446200419\\
31	69.8498569389295\\
32	69.8688733272201\\
33	69.8937670461917\\
34	69.8750186649387\\
35	69.7708663178535\\
36	69.8038126317793\\
37	69.743899825743\\
38	69.8827614684405\\
39	69.9403429554591\\
40	69.9798824456177\\
41	69.8735481500479\\
42	69.8656637793999\\
43	69.9051433514613\\
44	69.9444155982259\\
45	69.93792283883\\
46	69.8885473582753\\
47	69.8769474837554\\
48	69.8323897848251\\
49	69.8554771095653\\
50	69.9182853035033\\
51	69.9670967630742\\
52	69.9234775299385\\
53	69.9224757333274\\
54	69.9170970909433\\
55	69.9626748765597\\
56	69.9457781386838\\
57	69.9506502568368\\
58	69.9097230826976\\
59	69.8843058069243\\
60	69.8836228712756\\
61	69.894292866978\\
62	69.9421695844781\\
63	69.9174488988659\\
64	69.9366778030092\\
65	69.9127984111116\\
66	69.9178088451057\\
67	69.9290963739563\\
68	69.9357542064807\\
69	69.9331334081979\\
70	69.8977500058327\\
71	69.8873074712877\\
72	69.8620500468974\\
73	69.895960679078\\
74	69.9064788847296\\
75	69.9249556909015\\
76	69.92219666785\\
77	69.9133337060059\\
78	69.9206507190637\\
79	69.927710777536\\
80	69.9416052271402\\
81	69.9460850160328\\
82	69.926870243636\\
83	69.9042181161704\\
84	69.9161937797699\\
85	69.9256906151608\\
86	69.9408619677356\\
87	69.9322171545501\\
88	69.9498680451498\\
89	69.9335972343818\\
90	69.9509210097967\\
91	69.9530893955819\\
92	69.9625410549721\\
93	69.946391139173\\
94	69.9274605646506\\
95	69.9384385126711\\
96	69.9308391966178\\
97	69.9388771753745\\
98	69.9307093947062\\
99	69.9453690727629\\
100	69.9287857522208\\
};
\addlegendentry{mid}

 \addplot [color=green, line width=1.0pt]
  table[row sep=crcr]{%
0	69.958225249148\\
1	69.9582138500062\\
2	69.9582152678939\\
3	69.9582303014454\\
4	69.9582276232794\\
5	69.9582205657277\\
6	69.9582138482208\\
7	69.9582058795341\\
8	69.9581974615566\\
9	69.958186955549\\
10	69.9581822458707\\
11	69.958172990448\\
12	69.9581747289937\\
13	69.958180932453\\
14	69.9581985465699\\
15	69.9582028982883\\
16	69.9581959738672\\
17	69.9581996334946\\
18	69.9582033279668\\
19	69.9582136378374\\
20	69.9582126436981\\
21	69.8694646529975\\
22	69.8669657930684\\
23	69.8653359828856\\
24	69.8637284253578\\
25	69.8622593358755\\
26	69.8708475880516\\
27	69.9105248387492\\
28	69.9042916990526\\
29	69.9165912337141\\
30	69.9341874292862\\
31	69.931803053181\\
32	69.8821302014298\\
33	69.879325909773\\
34	69.8566886343529\\
35	69.8521089689086\\
36	69.8437468666595\\
37	69.834900427997\\
38	69.86763368687\\
39	69.8907191690104\\
40	69.9233956912644\\
41	69.9645845536582\\
42	69.9570419590228\\
43	69.933604784535\\
44	69.9521053104389\\
45	69.9387132803477\\
46	69.9533427260852\\
47	69.9215154272104\\
48	69.9121078429919\\
49	69.9241557340728\\
50	69.9218729150135\\
51	69.9233267590514\\
52	69.956643561776\\
53	69.9749863564565\\
54	69.9679600055056\\
55	69.9737478651936\\
56	69.9640529478505\\
57	69.9777419814152\\
58	69.9504339885136\\
59	69.9169961114249\\
60	69.9217846093197\\
61	69.9197153924262\\
62	69.9137564606056\\
63	69.9222780151806\\
64	69.9155400884578\\
65	69.9239096778457\\
66	69.9276428642641\\
67	69.9280037686606\\
68	69.9431704651052\\
69	69.9374521259427\\
70	69.918831076952\\
71	69.9111086134981\\
72	69.903573465311\\
73	69.893230565304\\
74	69.900222152047\\
75	69.8999272532962\\
76	69.9036311025877\\
77	69.9089214162645\\
78	69.9134175099355\\
79	69.9308464139818\\
80	69.935067228337\\
81	69.9194467268056\\
82	69.9292892649329\\
83	69.9280404622104\\
84	69.9269766357762\\
85	69.9210946691366\\
86	69.9144170378954\\
87	69.9180467631657\\
88	69.9185449510127\\
89	69.9203381608356\\
90	69.931821634207\\
91	69.9459387904694\\
92	69.9381866096047\\
93	69.9388881078266\\
94	69.9374498510127\\
95	69.9404156918558\\
96	69.9391270986157\\
97	69.9284378388371\\
98	69.9267264083292\\
99	69.9269950912186\\
100	69.9215892726112\\
};
\addlegendentry{end}

\addplot [color=black, line width=1.0pt]
  table[row sep=crcr]{%
0	69.9581200687398\\
1	69.9581319488398\\
2	69.9583376484376\\
3	69.9582759032365\\
4	69.9580910383862\\
5	69.9581974426274\\
6	69.958351254045\\
7	69.9582082537371\\
8	69.9580941699918\\
9	69.9582641517402\\
10	69.958332685827\\
11	69.9581471281189\\
12	69.9581275846246\\
13	69.9583170510495\\
14	69.9582873275805\\
15	69.9581059602232\\
16	69.9581825584828\\
17	69.9583447518986\\
18	69.9582263893681\\
19	69.9580932293788\\
20	69.957911266077\\
21	69.8699150102191\\
22	69.8667919406098\\
23	69.8652835851496\\
24	69.8548709096314\\
25	69.866195557401\\
26	69.8792616823452\\
27	69.892890492243\\
28	69.9025752783614\\
29	69.8994121844315\\
30	69.9050601588024\\
31	69.909942019302\\
32	69.8886128192685\\
33	69.8720545673935\\
34	69.8770992922536\\
35	69.8694633793686\\
36	69.8561705767806\\
37	69.8701530240461\\
38	69.891624736408\\
39	69.8992605963614\\
40	69.9127582275586\\
41	69.9329163390765\\
42	69.9388905702026\\
43	69.9390042927017\\
44	69.9460856395058\\
45	69.9472505193596\\
46	69.9395129874991\\
47	69.9375243307023\\
48	69.9420749578628\\
49	69.9432111018691\\
50	69.9400329959693\\
51	69.9425756993485\\
52	69.9542287478169\\
53	69.9577088479853\\
54	69.9484279565565\\
55	69.9505334189759\\
56	69.9602513048071\\
57	69.9489348790485\\
58	69.9324767337313\\
59	69.9374880389331\\
60	69.9402714023877\\
61	69.9239731887935\\
62	69.9159028383233\\
63	69.9252514658006\\
64	69.926950031021\\
65	69.9154578172904\\
66	69.9125322488695\\
67	69.9224777836226\\
68	69.9229529725791\\
69	69.9100611327352\\
70	69.9084889387798\\
71	69.9189053032798\\
72	69.9152683536265\\
73	69.9034208815401\\
74	69.9085630394311\\
75	69.918855364021\\
76	69.9137431318285\\
77	69.9088760853078\\
78	69.918258300734\\
79	69.9249061273324\\
80	69.920161722907\\
81	69.9193467052952\\
82	69.926879769838\\
83	69.9293369845988\\
84	69.9248921592446\\
85	69.9256273290627\\
86	69.9314080156503\\
87	69.9312269805568\\
88	69.9279205107878\\
89	69.9312270653714\\
90	69.9350557749903\\
91	69.9320550788727\\
92	69.9310567549304\\
93	69.9353635612508\\
94	69.9345570466457\\
95	69.9299234960856\\
96	69.9315786919712\\
97	69.9341623440187\\
98	69.9295230039705\\
99	69.9261385974443\\
100	69.9298570923973\\
};
\addlegendentry{new}

\end{axis}
\end{tikzpicture}%